\newcommand\webcite[1]{\texttt{\def~{\~{}}#1}}
\documentclass{article}
\usepackage{amsmath,amsfonts,amsthm,epsfig,color}
\newtheorem{theorem}{Theorem}
\newtheorem{lemma}[theorem]{Lemma}

\newtheorem{corollary}[theorem]{Corollary}
\theoremstyle{definition}
\newtheorem{remark}[theorem]{Remark}

\newcommand\bb[1]{\bigl(#1\bigr)}

\newcommand\dd{\,d}
\newcommand{\mucs}{$\mu$-continuity set}
\newcommand\sss{{\mathcal S}}
\newcommand\ddd{\partial}
\newcommand\vxs{\ensuremath{\mathcal V}}
\newcommand\xs{\ensuremath{{\bf x}_n}}
\newcommand\xss{\ensuremath{(\xs)_{n\ge 1}}}
\newcommand\xn{x^{(n)}}
\newcommand\gnkx{\gnxx{\kk}}
\newcommand\gnxx[1]{\ensuremath{G^\vxs(n,#1)}}
\newcommand\kk{\kappa}
\newcommand\pij{p_{ij}}
\newcommand\mui{\mu\{i\}}

\newcommand\Prk{{\mathop{\mathbb P_\kk}\nolimits}}

\newcommand\Prkx{{\mathop{\mathbb P_{\kk,x}}\nolimits}}
\newcommand\Pow{{\mathcal P}}
\renewcommand\Pr{{\mathop{\mathbb P{}}\nolimits}}
\newcommand\eps{\varepsilon}
\newcommand\downto{\searrow}
\newcommand\upto{\nearrow}
\newcommand\la{\lambda}
\newcommand\ka{\kappa}
\newcommand\pge[1]{\Psi_{\ge #1}}
\newcommand\pgem{\pge{k-1}}
\newcommand\pget{\pge{t}}
\newcommand\pgek{\pge{k}}
\newcommand\Prl{{\Pr_\la}}
\newcommand\br{{\mathcal B}}
\newcommand\ca{{\mathcal A}}
\newcommand\ps{{\mathcal P}}
\newcommand\cD{{\mathcal D}}
\newcommand\cT{{\mathcal T}}
\newcommand\cB{{\mathcal A}'}
\newcommand\cW{{\mathcal W}}
\newcommand\Ga{{\Gamma}}
\newcommand\tG{{\tilde G}}
\newcommand\whp{whp}
\newcommand\ev{{\mathcal E}}
\newcommand\el{{\mathcal L}}
\newcommand\rb{{\mathcal R}}
\newcommand\of{\circ}
\renewcommand\={:=}
\newcommand\Po{\operatorname{Po}}
\newcommand\E{\operatorname{\mathbb E}{}}
\newcommand\Es[1]{\operatorname{\mathbb E}_{#1}}
\newcommand\diff{\bigtriangleup}
\newcommand{\brs}{\brp_{\le 2L+1}}
\newcommand{\Ms}{{\mathcal M}_{\le 2L+1}}
\newcommand{\Mle}[1]{{\mathcal M}_{\le #1}}
\newcommand{\bp}{\beta^+}
\newcommand{\brp}{\br^+}
\newcommand{\vp}{{\bf p}}
\newcommand\pto{\overset{\mathrm{p}}{\to}}

\newcommand\Aeps{A^{(\eps)}}
\newcommand\kkeps{\kk^{(\eps)}}
\newcommand\gp{\gamma^+}
\newcommand\feps{f^{(\eps)}}
\newcommand\lac{\la_{\mathrm c}}
\newcommand\cc{{\mathrm c}}

\begin{document}
\title{The $k$-core and branching processes}
\date{February 4, 2007}

\author{Oliver Riordan%
\thanks{Royal Society Research Fellow, Department of Pure Mathematics
and Mathematical Statistics, University of Cambridge, Cambridge CB3 0WB}}
\maketitle

\begin{abstract}
The {\em $k$-core} of a graph $G$ is the maximal subgraph of $G$
having minimum degree at least $k$.
In 1996, Pittel, Spencer and Wormald found the threshold $\lac$ for the emergence
of a non-trivial $k$-core in the random graph $G(n,\la/n)$, and the asymptotic
size of the $k$-core above the threshold. We give a new proof of this result using
a local coupling of the graph to a suitable branching process. This proof extends to a general
model of inhomogeneous random graphs with independence between the edges.
As an example, we study the $k$-core in a certain power-law or `scale-free'
graph with a parameter $c$ controlling the overall density of
edges. For each $k\ge 3$, we find the threshold value of $c$ at which the $k$-core
emerges, and the fraction of vertices in the $k$-core when $c$ is $\eps$ above the threshold.
In contrast to $G(n,\la/n)$, this fraction tends to $0$ as $\eps\to 0$.
\end{abstract}

\section{Introduction}

The {\em $k$-core} of a graph $G$ is the maximal subgraph of $G$ with
minimum degree at least $k$; if $G$ has no such (non-trivial) subgraph,
then the $k$-core of $G$ is empty. This concept was introduced
by Bollob\'as~\cite{BBkcore}, in the context of finding large $k$-connected
subgraphs of random graphs. As edges are added one by one to a graph,
the $k$-core grows; in particular, it is empty up to some point, and then
non-empty (and often large). The question of when the $k$-core
emerges in a random graph (or random graph process) also arose
in the context of finding the chromatic number of sparse random graphs;
more specifically, let $G(n,p)$ be the Erd\H os-R\'enyi
random graph with $n$ vertices, in which the possible edges are present
independently, each with probability $p$. If a graph has no $k$-core,
it is $k$-colourable; Chv\'atal~\cite{Chv} used this to show
that $G(n,2.88/n)$ is \whp\ 3-colourable.
Here, as usual, an event holds  {\em with high probability},
or \whp,
if it holds with probability $1-o(1)$ as $n\to\infty$.

It is natural to ask: for
$k\ge 3$ fixed, what is
the critical value $\lac=\lac(k)$ of $\la$ above which a
(non-empty) $k$-core first appears \whp\ in $G(n,\la/n)$?

There is a natural
`guess' as to the answer, given by a branching process:
let $X_\la$ be a Galton--Watson branching process that starts with a
single particle $x_0$ in generation zero, where the number of
children of each particle has a Poisson distribution with mean $\la$,
and these numbers are independent for different particles.
Of course, $X_\la$ provides a good `local' model of the neighbourhood
of a vertex of $G(n,\la/n)$.
Let $\brp$ be the event that $x_0$ has at least $k$
children each of which has at least $k-1$ children each of which has at
least $k-1$ children each of which \dots, i.e., that $x_0$ is in a $k$-regular tree
contained in $X_\la$. (We reserve the notation $\br$ for an associated event
whose role in the analysis is more fundamental.)
Let $\bp(\la)$ be the probability that $X_\la$ has the property $\brp$.
One might expect that, up to probability $o(1)$, a vertex is in the $k$-core
if its neighbourhoods up to a suitable distance have a property
corresponding to $\brp$, and thus that the fraction
of vertices of $G(n,\la/n)$ in the $k$-core is $\bp(\la)+o(1)$;
this turns out to be the case.

Pittel, Spencer and Wormald~\cite{PSW} showed that,
except at the critical point, the number of vertices in the $k$-core
of $G(n,\la/n)$ is indeed $\bp(\la)n+o_p(n)$. In particular,
the threshold $\lac$ for the emergence of the $k$-core is
$\lac=\inf\{\la:\bp(\la)>0\}$.
Recently, simpler proofs of this result have been given, as well
as generalizations to various other contexts; see, for example,
\cite{Coop,Molloy2005,FR,CW,JL,DarlingNorris}.
As far as we are aware, although the branching process heuristic
was described already in \cite{PSW},
none of these proofs works by directly coupling the neighbourhoods of a vertex
of the graph with the branching process; there seems to be a problem
relating the inescapably global property of lying in the $k$-core
to a simple local property.
Here we give a new proof that does proceed in this way,
using a carefully chosen (and not very simple) local property.

As the proof will be a little involved, we give a very rough outline
to motivate what follows: we shall define a certain event
$\ca$ depending on the first $D$ generations of $X_\la$, where $D=o(\log
n)$ but $D/\log\log n\to\infty$, such that $\Pr(\ca)$ is almost as large
as $\Pr(\brp)$, and $\Pr(\brp\mid \ca)=1-o(n^{-1})$. As $\ca$ is a `local'
event, the number of vertices in $G(n,\la/n)$ with (neighbourhoods
having the equivalent of) property $\ca$ will be
$\Pr(\ca)n+o_p(n)$. Defining $\ca$ in the right way, we can show that,
given that a vertex has property $\ca$, with probability $1-o(n^{-1})$ it is the
root of a $k$-regular tree of height at most $D$ all of whose leaves have
property $\ca$.  But then \whp\ {\em every} vertex with property $\ca$ is the root of
such a tree, and the union of these trees has minimum degree at least
$k$. Most of the work will be in the analysis of the branching
process: having found the right event, the translation to the graph
will be relatively straightforward.

Unfortunately, the event $\ca$ we have to work with is rather
complicated.  Roughly speaking, we look for a tree $T$ of height $D$ that
branches a little faster than a $k$-regular tree. Having found such a tree,
if we select leaves of $T$ independently
with an appropriate probability (around $\beta(\la)$),
it is very likely that we can find a $k$-regular tree whose
leaves are a subset of the selected leaves. When
we explore the neighbourhood of a vertex $v$ in the graph, finding
a tree $T$ as above, we then wish to explore from each leaf $w$ of $T$
to see whether the neighbourhoods of $w$ have the property $\ca$.
In doing this, it turns out that we cannot afford to ignore the
edge along which we reached $w$ when exploring from $v$, and that
we must re-use parts of the tree $T$ in establishing that $w$
has property $\ca$. This, together with the need to achieve (approximate)
independence of the explorations from different leaves $w$,
places rather subtle constraints on the
properties $\ca$ we can use.

The next section is devoted to the study of the $k$-core in $G(n,\la/n)$.
In Section \ref{sec_nonunif}, we turn to inhomogeneous random graphs:
an advantage of the direct branching process approach is that
it extends easily to other random graph models in which
the neighbourhoods of a vertex can be modelled by a suitable
branching process. This includes the general sparse inhomogeneous model
of Bollob\'as, Janson and Riordan~\cite{kernels}.
Indeed, the present work started at the conference
Random Structures and Algorithms 2005 in Poznan,
when Alan Frieze asked whether we knew the size of the $k$-core in this model.
As a special case, we study the graph on $[n]$ in which edges
are present independently, and the probability of an edge between
$i$ and $j$ is $c/\sqrt{ij}$ for a parameter $c>0$.
We find the threshold (in terms of $c$) at which the $k$-core appears
in this graph, and the size of the $k$-core above this threshold; in contrast
to $G(n,\la/n)$, just above the threshold the $k$-core is small.

\section{The uniform case}\label{sec_unif}

We start with some standard basic preliminaries.
Let us write $\Prl$ for the probability measure associated to
the branching process $X_\la$; when discussing $X_\la$,
we shall use the terms `event' and `property of the branching
process' interchangeably, to refer to a (measurable) subset of the set
of all rooted,
unlabelled trees. When we say that a particle $x$ of the
branching process has a certain property, we mean that the process
consisting of $x$ (as the new root) and its descendants has the
property. Throughout we write $x_0$ for the initial particle
of $X_\la$. When $\la$ is clear from the context (or specified in the
notation for our probability measure), we may write $X$ for $X_\la$.
Many of the events and quantities we shall define below depend on $k$.
Almost always, we regard $k\ge 2$ as fixed and suppress this dependence.

Let $\pget(\la)=\Pr(\Po(\la)\ge t)$ denote the probability that a Poisson random variable with
mean $\la$ takes a value that is at least $t$.
Let $\br_d$ be the event that $X$ contains a $(k-1)$-ary tree of height
$d$ with $x_0$ as the root, and let $\br=\lim_{d\to\infty}\br_d$
be the event that $X$ contains an infinite $(k-1)$-ary tree
with root $x_0$.
Then $\Prl(\br_0)=1$. Also,
each particle in the first generation of $X_\la$ has probability
$\Prl(\br_d)$ of having property $\br_d$. As these events are independent
for different particles, the number of particles in the first generation
with property $\br_d$ has a Poisson distribution with mean
$\la\Prl(\br_d)$. Thus, $\Prl(\br_{d+1})=\pge{k-1}\bb{\la\Prl(\br_d)}$.

Since $\pge{k-1}(\la x)$ is a continuous, increasing function of $x$, it follows
that $\Prl(\br)=\lim_{d\to\infty}\Prl(\br_d)$ is given by the maximum solution
$p$ to the equation $p=\pge{k-1}(\la p)$. We denote this solution by $\beta(\la)$.
As in the introduction, let $\brp$ be the event that $X$ contains an infinite $k$-regular
tree with root $x_0$, i.e., that $x_0$ has at least $k$ children with property
$\br$; we shall write $\brp_d$ for the corresponding event depending
only on the first $d$ generations, i.e., that $x_0$ has at least $k$ children
with property $\br_{d-1}$. Note that
\[
 \bp(\la)\=\Prl(\brp)=\pge{k}\bb{\la\Prl(\br)}=\pgek\bb{\la \beta(\la)}.
\]
Let $\lac=\inf\{\la:\beta(\la)>0\}$. It is easy to see that (in this uniform case)
the functions $\beta(\la)$ and $\bp(\la)$ have a jump at $\la=\lac$ (for $k\ge 3$),
and are continuous and strictly increasing for $\la\ge\lac$.

Pittel, Spencer and Wormald~\cite{PSW} proved the following result,
with sharper error estimates. As usual, if $A_n$ is a sequence of random variables,
we write $A_n=o_p(f(n))$ if $A_n/f(n)\to 0$ in probability, i.e.,
if $|A_n|\le \eps f(n)$ \whp\ for any $\eps>0$.

\begin{theorem}[\cite{PSW}]\label{t1}
Let $k\ge 3$ be fixed, and define $\lac$ and $\bp(\la)$ as above. If $\la\ne\lac$
is constant, then the number of vertices of $G(n,\la/n)$ in the $k$-core
is $\bp(\la)n+o_p(n)$ as $n\to\infty$.
\end{theorem}

Note that the $2$-core of a graph $G$ is just the union of the cycles in $G$,
together with all paths joining two cycles.
It is easy to see that a large $2$-core first appears in $G(n,\la/n)$
at the same time as the giant component, i.e., when $\la>\lac(2)=1$, although there may be
a small $2$-core when $\la<1$. With the weak error bounds above, Theorem \ref{t1}
holds for $k=2$ as well. For us, there will be no difference between
the cases $k=2$ and $k\ge 3$, although the latter is the more interesting.

The proof of the upper bound is easy; we postpone this to the end of
the section.  Our proof of the lower bound will require considerable
preparation. Although the functions $\beta(\la)$ and $\bp(\la)$ are
continuous except at $\la=\lac$, we shall avoid using this fact, with
an eye to generalizations.

If $\ev_1$, $\ev_2$ are properties of the branching process, with $\ev_1$ depending only on the
first $d$ generations, let $\ev_1\of\ev_2$ denote the event that $\ev_1$ holds
if we delete from $X$ all particles in generation $d$ that do not have property $\ev_2$.
For example, with $\br_1$ the property of having at least $k-1$ children, as above,
$\br_1\of\br_1=\br_2$, the property of having at least $k-1$ children with at least
$k-1$ children (and perhaps other children with fewer than $k-1$ children).

Let $\rb_d$ be the event that $\br_d$ holds in a {\em robust} manner, meaning
that $\br_d$ holds even after any single particle in generation $d$ is deleted.
Then, as $\br=\br_d\of\br$, the event $\rb_d\of\br$ is the event that $\br$ holds 
after any single particle in generation $d$ is deleted. Of course,
whenever a particle is deleted,
so are its descendants.

\begin{lemma}\label{l1}
If $\la\mapsto\beta(\la)$ is continuous at $\la$ then
\[
 \Prl(\rb_d\of\br)\upto \beta(\la)\=\Prl(\br)
\]
as $d\to\infty$.
\end{lemma}
\begin{proof}
Fix $0<\eps<1$. Consider the natural coupling of $X_\la$ and $X_{(1-\eps)\la}$, obtained
by constructing $X_\la$, and then deleting each edge (of the rooted tree)
independently with probability $\eps$, and taking for $X_{(1-\eps)\la}$ the set of particles
still connected to the root. In this coupling, if $\br$ does not hold for $X_\la$,
it certainly does not hold for $X_{(1-\eps)\la}$. Furthermore, if $\br\setminus (\rb_d\of\br)$
holds for $X_\la$, then there is some particle $x$ in generation $d$ such that
if $x$ is deleted, then $\br$ no longer holds. But the probability that $x$ is
not deleted when passing to $X_{(1-\eps)\la}$ is exactly $(1-\eps)^d$.
It follows that
\[
 \Pr_{(1-\eps)\la}(\br) \le \Prl(\rb_d\of\br)+(1-\eps)^d.
\]
Since $\rb_d\subset \rb_{d+1}$, the sequence $\Prl(\rb_d\of\br)$ is increasing. Taking the limit
of the inequality above,
\[
 \beta\bb{(1-\eps)\la} = \Pr_{(1-\eps)\la}(\br) \le \lim_{d\to\infty} \Prl(\rb_d\of\br).
\]
Letting $\eps\to 0$, the lemma follows.
\end{proof}

It will often be convenient to {\em mark} some subset of the particles in generation $d$.
If $\ev$ is an event depending on the first $d$ generations, then we write
$\ev\of M$ for the event that $\ev$ holds after deleting all unmarked particles in generation
$d$. We write $\Pr^p_\la(\ev\of M)$ for the probability that $\ev\of M$ holds
when, given $X_\la$, we mark the particles in generation $d$ independently with probability $p$.
We suppress $d$ from the notation, since it will be clear from the event $\ev$.
Let
\[
 r(\la,d,p) = \Pr^p_\la(\rb_d\of M).
\]

\begin{lemma}\label{l2}
Let $\la_1< \la_2$ be fixed, with $\la\mapsto \beta(\la)$ continuous at $\la_1$. Then there is
a $d$ such that
\begin{equation}\label{re}
 r\bb{\la_2,d,\beta(\la_1)} \ge \beta(\la_1).
\end{equation}
\end{lemma}
\begin{proof}
Let us construct a branching process $X'$ as follows. Start with $X_{\la_1}$. Then, independent
of $X_{\la_1}$, add a $\Po(\la_2-\la_1)$ number of `extra' children of the initial particle. Each extra child
then has descendants as in $X_{\la_1}$. 
Clearly, we may consider $X'$ as a subset of $X_{\la_2}$.
Given $d\ge 1$, let us mark each particle in generation $d$ if it has property $\br$.
Now the descendants of a particle in generation $d$ have the distribution of $X_{\la_1}$,
with independence for different particles. Thus, given the first $d$ generations of $X'$,
we mark each particle in generation $d$ independently with probability $p=\beta(\la_1)$.
Hence, $r(\la_2,d,p)$ is at least the probability that $\rb_d\of \br$ holds in $X'$.

Each extra child has property $\br$ independently with probability $p$.
Hence, the event $\ev_1$ that there are at least $k$ extra children with property $\br$
has probability $\delta=\pgek((\la_2-\la_1)p)>0$.
Now $\ev_1$ is independent of $X_{\la_1}$.
Also, if $\ev_1$ holds in $X'$, then so does
$\rb_1\of\br$, and hence $\rb_d\of\br$ for any $d$.
Thus,
the probability that $\rb_d\of\br$ holds in $X'$ is at least
\[
 r_d=1- (1-\delta)\bb{1-\Pr_{\la_1}(\rb_d\of\br)}.
\]
By Lemma~\ref{l1}, as $d\to\infty$ we have
$\Pr_{\la_1}(\rb_d\of\br)\to \beta(\la_1)=p$, so
\[
 r_d\to 1-(1-\delta)(1-p)>p,
\]
and there is a $d$ with $r_d\ge p$, completing the proof.
\end{proof}

We shall use robustly branching trees, i.e., rooted trees with height
$d$ having property $\rb_d$, in two ways. The first is the obvious way:
we shall use the fact that at least two leaves must be deleted in order
to destroy all branching subtrees, i.e., all subtrees with property
$\br_d$. The second is less direct, and is described in the next lemma,
for which we first need a definition.

We say that a finite rooted tree has the property $\Mle{d}$
if all its leaves are at distances between $1$ and $d$ from the root, and every non-leaf
has degree at least~$k$.

\begin{lemma}\label{brp}
Let $\cT$ be a rooted tree of height $d+1$ that is minimal with respect
to having property $\rb_d\of\br_1$.
Let $w$ be a vertex of $\cT$ at distance $d$ from the root, and let $y$
be the parent of $w$.
If the graph $\cT$ is regarded as a rooted tree with root $w$,
then it has a subtree $\cW$ with the property $\Mle{2d+1}$,
such that all leaves of $\cW$ are leaves of $\cT$.
\end{lemma}

\begin{figure}[htb]
\centering
\input{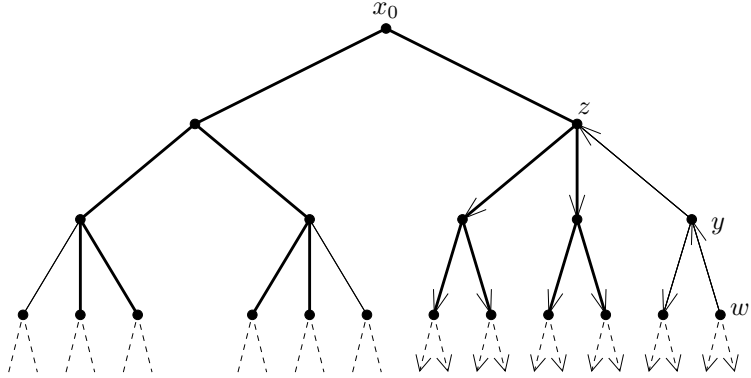}
\caption{A tree $\cT$ with property $\rb_3\of\br_1$, where $k=3$. The last generation is drawn with
dotted lines; the solid tree $\cT'$ has property $\rb_3$. The thick lines are a minimal subtree
$\cT''$ of $\cT'-w$ with property $\br_3$. The arrows consist of the shortest path $wz$ from
$w$ to $\cT''$ and all descendants of vertices in this path. Taking $w$ as the root,
they form a tree with property $\Mle{5}$.} \label{fig1}
\end{figure}

\begin{proof}
The proof is illustrated in Figure \ref{fig1}. Although it is not clear that
a written proof adds anything to the figure, to be formal we include one.

Let $\cT'$ denote the first $d$ generations
of $\cT$, a tree with property $\rb_d$ in which (by minimality)
every non-leaf has at least $k-1$ children. By definition of the property $\rb_d$,
deleting $w$ from $\cT'$ leaves a tree with property $\br_d$; let
$\cT''$ be a minimal subtree of $\cT'-w$ with property $\br_d$, so $\cT''$ is
a $(k-1)$-ary tree.
In $\cT$, there is a unique
path $P$ from $w$ to $\cT''$, meeting $\cT''$ at $z$, say.
Consider the tree $\cW$ obtained by taking $P$
together with the descendants in $\cT$ of every vertex of $P$. Regarding
this tree as rooted at $w$, the children
of $z$ in $\cW$ are exactly its children in $\cT''$, of which
there are $k-1$. Furthermore, any vertex $a$ of $P$ other than $z$ and $w$
has the same number of children in $\cW$ as in $\cT$ (we lose one neighbour
on $P$ and gain the other), which is at least $k-1$. All other vertices
of $\cW$ have the same children in $\cW$ as in $\cT$, except
for $w$, which has one additional child $y$ in $\cW$. Thus
all non-leaves of $\cW$ have at least $k-1$ children, the root $w$ has
at least $k$, and all leaves of $\cW$ are leaves of $\cT$.
Observing finally that the height of $\cW$ is at most $2d+1$, the result follows.
\end{proof}

Let us briefly recall Harris' Lemma. Let $X$ be a finite set (here
consisting of some or all possible edges of $G(n,p)$),
and let $X_p$ be a random subset of $X$ formed by selecting each
element independently with probability $p$. An event
$\ca\subset\Pow(X)$ is {\em increasing} if $A\in \ca$ and $A\subset B\subset
X$ imply $B\in \ca$, and {\em decreasing} if $A\in \ca$ and $B\subset A$
imply $B\in \ca$.
\begin{lemma}
If $\ca_1,\ca_2\subset \Pow(X)$ are increasing events and $0<p<1$, then
\[
 \Pr(X_p\in \ca_1\cap \ca_2) \ge \Pr(X_p\in \ca_1)\Pr(X_p\in \ca_2).
\]
\end{lemma}
In other words, increasing events are positively correlated. Of course,
it follows that decreasing events are also positively correlated.
This was proved by Harris in 1960~\cite{Harris}, and rediscovered
by Kleitman~\cite{Kleitman}.

We shall need one simple lemma concerning $G(n,p)$; this is because we wish
to prove a result about the neighbourhoods of {\em every} vertex, and cycles
within these neighbourhoods will cause problems. To avoid this, we consider
the random graph $\tG=\tG(n,\ell,\la)$ whose distribution is that of
$G(n,\la/n)$ conditioned on the absence of any cycles of length at most $\ell$.
Although globally $G(n,\la/n)$ is very likely to contain cycles, locally,
conditioning on their absence makes little difference: after this conditioning, given that
a certain not too large set of edges is present, and that a certain other
set of edges is absent, the probability that another edge is present
is close to $\la/n$, as long as this edge would not complete a `known' cycle.

\begin{lemma}\label{l_small}
Let a constant $\la>0$ and a function $\ell=\ell(n)=o(\log n)$ be given.
If $n$ is large enough then, whenever $E_0$, $E_1$ and $\{e\}$ are disjoint
sets of possible edges of $\tG=\tG(n,\ell,\la)$ with $|E_1|\le n^{1/3}$ such
that $E_1\cup \{e\}$ contains no cycle of length at most $\ell$, we have
\[
 (1-n^{-1/4})\la/n \le \Pr\bb{e\in E(\tG) \mid E_1\subset E(\tG) \subset E_0^\cc} \le \la/n.
\]
\end{lemma}
\begin{proof}
Let $G'$ denote the (distribution of the) random graph $\tG$ conditioned
on the presence of every edge in $E_1$ and the absence of every edge in $E_0$,
so the probability we wish to bound is $\Pr(e\in G')$.
Note that $G'$ may be described as follows: start from $G(n,\la/n)$, and first
condition on the presence of the edges
in $E_1$ and the absence of the edges in $E_0$. At this point,
each edge of $E_2=(E_0\cup E_1)^\cc$
is present independently with probability $\la/n$. Writing $\ps(S)$ for
the power-set of a set $S$, we next condition on the
event $\cD\subset\ps(E_2)$ that the set of edges $E$ from $E_2$ present is such that
the whole graph does not contain a short (length at most $\ell$) cycle,
i.e., that $E\cup E_1$ does not contain a short cycle.
In terms of the edges in $E$, 
the event $\cD$ is a decreasing event.
The upper bound (which we shall not in fact use) is thus  essentially trivial:
as the event $e\in E$ is increasing, by Harris' Lemma
\[
 \Pr(e\in G')=\Pr(e\in E\mid \cD)\le \Pr(e\in E)=\la/n.
\]

For the lower bound, let us call a set $P\subset (E_0\cup E_1\cup\{e\})^\cc$
a {\em pre-cycle} if it is minimal subject to $P\cup E_1\cup \{e\}$ containing
a short cycle that includes the edge $e$. As legal configurations in the graph
$G'$ have the same relative probabilities in $G'$ as in $G(n,\la/n)$, the conditional
probability
$\Pr(e\in G'\mid G'\setminus\{e\})$ is $0$ if $G'\setminus\{e\}$ contains
a pre-cycle, and $\la/n$ otherwise.
Thus,
\begin{eqnarray*}
 \Pr(e\in G') &=& 
 \E\bb{\Pr(e\in G'\mid G'\setminus\{e\})} \\
 &=&  (\la/n)\Pr(G'\setminus\{e\}\hbox{ contains no pre-cycle}) \\
 &=& (\la/n)\Pr(E\hbox{ contains no pre-cycle}\mid \cD) \\
 &\ge& (\la/n)\Pr(E\hbox{ contains no pre-cycle}),
\end{eqnarray*}
where the last step is again from Harris' Lemma.
Thus, it suffices to show that with probability at least
$1-n^{-1/4}$ the set $E$ contains no pre-cycle. But
any pre-cycle must include a path of length $1\le t\le \ell$ joining
two vertices in $V$, the set of endpoints of edges in $E_1\cup\{e\}$.
Since $|V|\le 2n^{1/3}+2$, the expected number
of such paths is at most
\[
 \sum_{t=1}^\ell  \binom{|V|}{2} n^{t-1}(\la/n)^t
 \le \frac{4n^{2/3}}{n} \sum_{t=1}^\ell \la^t \le 4n^{-1/3}O(1)^{o(\log n)} \le n^{-1/4}
\]
if $n$ is large enough. Hence, the probability that such a path is present
is at most $n^{-1/4}$, and so is the probability that $E$ contains a pre-cycle,
completing the proof.
\end{proof}

We are now ready to prove Theorem~\ref{t1}.

\begin{proof}[Proof of Theorem~\ref{t1}]
We start with the more difficult bound, the lower
bound on the size of the $k$-core.
Note that $\beta(\la)$ and $\bp(\la)$ are continuous except at $\la=\lac$. (We
shall not use this fact elsewhere in this proof.)
Hence it suffices to prove that if $\la<\la'$, then \whp\ the $k$-core
of $G(n,\la'/n)$ contains at least $\bp(\la)n$ vertices. Let us fix $\la<\la'$.
We may assume that $\la\ge \lac$, as otherwise there is nothing to prove.
Fix any $\la<\la_2<\la'$ with $\beta(\cdot)$ continuous at $\la_2$. (An increasing
function has at most countably many discontinuities.)
Letting $\la_1\upto \la_2$, we have $\beta(\la_1)\upto \beta(\la_2)$.
Hence,
\[
 \pgem\bb{\la'\beta(\la_1)}\upto\pgem\bb{\la'\beta(\la_2)} > \pgem\bb{\la_2\beta(\la_2)} = \beta(\la_2).
\]
Thus we may choose $\la<\la_1<\la_2$ so that
\begin{equation}\label{close}
 \pgem\bb{\la'\beta(\la_1)} > \beta(\la_2).
\end{equation}
We may also assume that $\beta(\la)$ is continuous at $\la_1$.

By Lemma \ref{l2} there is a constant $d$ such that \eqref{re} holds; fix
such a $d$. Once $\beta(\la)$ is positive, it is strictly
increasing, so $\beta(\la_2)>\beta(\la_1)$.
Hence, by Lemma \ref{l1},
there is a $d_1$ such that $\Pr_{\la_2}(\rb_{d_1}\of\br)>\beta(\la_1)$.

Set
\[
 \eta=\min\big\{\Pr_{\la_2}(\rb_{d_1}\of\br)-\beta(\la_1),k^{-3d}\big\}>0.
\]
As $\rb_{d_1}\of\br$ is measurable, there is an integer $L$ and an event $\el_1$ depending
only on the first $L$ generations of the branching process such that
$\Pr_{\la_2}\bb{(\rb_{d_1}\of\br)\diff \el_1}\le \eta^2/2$,
where $\diff$ denotes symmetric difference.
We may and shall assume that $L>d_1$.
Writing $1_\ev$ for the indicator function of an event $\ev$, and $X[L]$
for the first $L$ generations of the branching process, we have
\begin{equation}\label{ebd}
 \eta^2/2 \ge \Pr_{\la_2}\bb{(\rb_{d_1}\of\br)^\cc \cap \el_1}
 = \Es{\la_2}\Bigl( 1_{\el_1} \Pr_{\la_2}\bb{(\rb_{d_1}\of\br)^\cc \mid X[L]} \Bigr),
\end{equation}
where $\Es{\la_2}$ is the expectation corresponding to $\Pr_{\la_2}$.
Set
\[
 \el = \el_1 \cap \Bigl\{ \Pr_{\la_2}\bb{\rb_{d_1}\of\br \mid X[L]} \ge 1-\eta \Bigr\},
\]
noting that the event $\el$ depends only on the first $L$ generations of $X$.
Since the expectation appearing in \eqref{ebd} is at least $\eta\Pr_{\la_2}\bb{\el_1\setminus \el}$,
we have $\Pr_{\la_2}\bb{\el_1\setminus\el}\le \eta/2$, so
\begin{equation}\label{pl}
 \Pr_{\la_2}(\el) \ge \Pr_{\la_2}(\el_1)-\eta/2 \ge
 \Pr_{\la_2}(\rb_{d_1}\of\br)-\eta^2/2-\eta/2 \ge \beta(\la_1).
\end{equation}
Also, recalling that $L>d_1$, so $\rb_L\of\br\supset \rb_{d_1}\of\br$, whenever
$\el$ holds we have
\begin{equation}\label{cpl}
 \Pr_{\la_2}\bb{\rb_L\of\br\mid X[L]} \ge
 \Pr_{\la_2}\bb{\rb_{d_1}\of\br\mid X[L]}\ge 1-\eta\ge 1-k^{-3d}.
\end{equation}

Let $\ca_0=\el$, and for $t\ge 1$ set $\ca_t=\rb_d\of\ca_{t-1}$.
Thus, $\ca_t$ is a `recursively robust' version of the event $\br_{dt}\of\el$.
From the independence of the descendants of different particles in generation
$d$ of $X_{\la_2}$ we have
\[
 \Pr_{\la_2}(\ca_t) = r\bb{\la_2,d,\Pr_{\la_2}(\ca_{t-1})}.
\]
Thus, by \eqref{pl}, \eqref{re}, and induction on $t$, we have
\begin{equation}\label{a1}
 \Pr_{\la_2}(\ca_t)\ge \beta(\la_1)
\end{equation}
for every $t$.
Recalling that $\la_2<\la'$, let us note for later that
\begin{multline}\label{a2}
 \Pr_{\la'}(\br_1\of\ca_t) = \pgem\bb{\la'\Pr_{\la'}(\ca_t)}
 \ge \pgem\bb{\la'\Pr_{\la_2}(\ca_t)} \\ \ge \pgem\bb{\la'\beta(\la_1)} > \beta(\la_2),
\end{multline}
where the final inequality is from \eqref{close}.

If $\ev$ depends on the first $d$ generations of $X$, let $\Pr^p(\ev\of M\mid X[d])$
denote the conditional probability given the first $d$ generations that
$\ev\of M$ holds if we mark particles in generation $d$ independently with probability $p$.
Note that this conditional probability does not depend on $\la$.
From \eqref{cpl}, whenever $\ca_0=\el$ holds we have
\begin{equation}\label{c0}
 \Pr^{\beta(\la_2)}(\rb_L\of M\mid X[L]) = \Pr_{\la_2}(\rb_L\of\br\mid X[L])
 \ge 1-k^{-3d}.
\end{equation}
Let $\ev_t=\br_{dt}\of\rb_L$, so $\ev_0=\rb_L$ and $\ev_t=\br_d\of\ev_{t-1}$ for $t>0$.
We claim that for any $t\ge 0$, whenever $\ca_t$ holds, then
\begin{equation}\label{c1}
 \Pr^{\beta(\la_2)}(\ev_t\of M \mid X[dt+L]) \ge 1-k^{-(2^t+2)d}.
\end{equation}
For $t=0$, this is just \eqref{c0}. We prove the claim by induction on $t$; let
us assume that $t>0$, and that \eqref{c1} holds with $t$ replaced by $t-1$. We condition
on the first $dt+L$ generations of $X$, and assume that $\ca_t$ holds.
Since $\ca_t=\rb_d\of\ca_{t-1}$, there is a set $Y$ of particles in generation $d$ such that $\ca_{t-1}$
holds for each $y\in Y$, and $\rb_d\of M$ holds if we mark only the particles in $Y$.
As any tree witnessing $\rb_d$ contains a subtree witnessing $\rb_d$ in which each particle
has at most $k$ children, we may assume that $|Y|\le k^d$.
Now let us mark each particle in generation $dt+L$ independently with probability
$\beta(\la_2)$, and let $Y'$ be the set of particles $y$ in $Y$ such that
$\ev_{t-1}\of M$ holds for $y$. Each $y\in Y$ is included independently
in $Y'$, and, by the induction hypothesis, the probability that a particular $y\in Y$ is
included is at least $1-k^{-(2^{t-1}+2)d}$. Thus
\[
 \Pr(|Y\setminus Y'|\ge 2) \le \binom{|Y|}{2} \left(k^{-(2^{t-1}+2)d}\right)^2 \le k^{2d}k^{-(2^t+4)d}
 = k^{-(2^t+2)d}.
\]
From the definition of $\rb_d$, whenever $|Y\setminus Y'|\le 1$
then $\br_d$ holds if we keep in generation
$d$ only the particles in $Y'$. But then $\br_d\of\ev_{t-1}\of M=\ev_t\of M$
holds, proving \eqref{c1}. This proves the claim by induction.

\medskip
Given an event $\ev$, let $\brs\of \ev$ be the event that
there is a set $S$ of {\em targets} in generations between $1$ and $2L+1$
such that $\ev$ holds for each target, and in the tree consisting of the
targets and their ancestors every non-leaf (i.e., every non-target) has degree at least $k$;
in the notation of Lemma~\ref{brp}, we require the targets and their ancestors to form
a tree with property $\Mle{2L+1}$.
Note that we need not (and in general must not) take all
particles in generations between $1$ and $2L+1$ that have property $\ev$ to be targets.
Thus, if $\ev$ is increasing, so is $\brs\of\ev$.

Let $T=T(n)$
be any function with $T=o(\log n)$ and $T/\log\log n\to \infty$.
Note that, until now, the branching process events we have defined
do not depend on $n$.
Set $\ca=\brs\of \ca_T$.
Since $\brp_1\of\ca_T\subset \ca$, from \eqref{a1} we have
\[
 \Pr_{\la_2}(\ca)\ge \pgek\bb{\la_2\Pr_{\la_2}(\ca_T)} \ge \pgek\bb{\la_2\beta(\la_1)}
 > \pgek\bb{\la_1\beta(\la_1)} = \bp(\la_1).
\]

Let $s=2L+1+dT+L=o(\log n)$, so $\ca$ depends on the first $s$ generations of
the branching process. By standard properties of branching processes,
for any constant $\la$ we have
\begin{equation}\label{small}
 \Pr_{\la}\bb{|X[2s]|\ge n^{1/10}} = o(n^{-1}),
\end{equation}
where $|X[2s]|$ is the total size of the first $2s$ generations
of the branching process. Indeed, a simple calculation shows that for any $m$
we have $\Es{\la}(|X[t]|^m)=O(\la^{mt})$ as $t\to\infty$. Thus
$\Es{\la}\bb{|X[2s]|^{20}}=O\bb{\la^{40s}}=o(n)$, and \eqref{small} follows
by Markov's inequality.

With this branching-process preparation behind us, we are now ready to turn
to the graph $G(n,\la'/n)$.

Let $\tG=\tG(n,4s,\la')$ be the graph obtained from $G(n,\la'/n)$
by conditioning on the absence of any cycles of length at most $4s$.
Let $c_k(G)$ denote the number of vertices in the $k$-core of a graph $G$.
As $\{c_k(G)\ge x\}$ is an increasing event, while the absence of
short cycles is a decreasing event, we have
\[
 \Pr(c_k(\tG)\ge x) \le \Pr\bb{c_k(G(n,\la'/n))\ge x}
\]
for any $x$. Hence it suffices to prove a lower bound for $c_k(\tG)$
that holds \whp\ instead of the corresponding bound for
$c_k(G)$. We shall work entirely with $\tG$.

Let $v$ be a random vertex of $\tG$, and explore its successive neighbourhoods
$\Ga_t(v)$, $t\le 2s$, in the usual way, where $\Ga_t(v)$ is the set of vertices of $\tG$
at graph distance $t$ from $v$. We start with $v$ `active' and all other vertices
`untested'. At each step we pick an `active' vertex $w$ closest to $v$, and test edges
from $w$ to untested vertices one by one, marking any neighbour of $w$ found in this way
as `active'. After testing all possible edges from $w$ we mark $w$ as `tested'.
By Lemma~\ref{l_small}, as long as we have reached at most $n^{1/3}$ vertices,
conditional on everything so far each test succeeds with probability
$(1+O(n^{-1/4}))\la'/n$; we never attempt
to test an edge that might complete a cycle. As the number of untested vertices is $n-O(n^{1/3})$,
we may couple the number of new neighbours of $w$ found with a Poisson distribution
with mean $\la'$ so that the two numbers agree with probability $1-O(n^{-1/4})$.

Let $G_v[t]$ be the subgraph of $\tG$ formed by the vertices within distance
$t$ of $v$, noting that for $t\le 2s$ this graph is by definition a tree.
Since, \whp, $X_{\la'}[2s]$ contains at most $n^{1/10}$ particles,
\whp\ $G_v[2s]$ contains at most $n^{1/10}$ vertices. It follows
from the coupling above that we may couple $G_v[2s]$ and $X_{\la'}[2s]$ so
as to agree in the natural sense \whp. Let us say
that $v$ has property $\ca$ if $G_v[2s]$ has property $\ca$ when viewed
as a branching process, an event that depends only on $G_v[s]$. We have shown that $v$ has property $\ca$
with probability $\Pr_{\la'}(\ca)+o(1)$.

To establish concentration of the number of vertices $v$ with property $\ca$,
we use a simple trick also used in \cite{kernels}.
Let $v$ and $w$ be independently chosen random vertices of $\tG$. Since the probability
that $v$ and $w$ are within distance $2s$ is $o(1)$, we may couple $G_v[s]$ and $G_w[s]$
with independent copies of $X_{\la'}[s]$ so as to agree \whp.
Hence, the probability that both $v$ and $w$ have property $\ca$
is $\Pr_{\la'}(\ca)^2+o(1)$. Writing $N$ for the number of vertices with property $\ca$, we thus
have
\[
 \E(N)=\Pr_{\la'}(\ca)n+o(n), \quad \E(N^2)=\Pr_{\la'}(\ca)^2n^2+o(n^2).
\]
Hence (by Chebyshev's inequality), $N/n$ converges in probability
to
\[
 \Pr_{\la'}(\ca) \ge \Pr_{\la_2}(\ca) \ge \bp(\la_1) > \bp(\la).
\]
Thus,
\[
 \Pr(N\le \bp(\la)n) =o(1).
\]
To complete the proof of the lower bound on the size of the $k$-core,
it thus suffices to show that $c_k(\tG)\ge N$ \whp. We do this by showing that
\whp\ {\em every} vertex with property $\ca$ is in the $k$-core.

Let $v$ be a random vertex of $\tG$. Let us test whether $v$ has property $\ca$
by exploring its neighbourhoods as follows. Throughout we regard the graph $G_v[2s]$
as a branching process with $v$ as initial particle; this is valid since the graph
contains no cycles. Thus we shall speak of the children and descendants of a vertex.
If $w$ is a vertex in generation $t$ (i.e., at graph distance $t$ from $v$),
we shall write $D(w)$ for the set of descendants of $w$ in generation $t+(dT+L)$.
Let $S_1$ and $S_2$ be initially empty sets of vertices. At each stage,
every vertex in $S_1$ will have property $\ca_T$ in the tree/branching
process rooted at $v$. In other words, $w$ and its descendants
will form a tree with property $\ca_T$. Recall that $\ca=\brs\of\ca_T$.

First, let us reveal the descendants of $v$ unto generation $dT+L+1$.
For each child $w$ of $v$, examine its descendants to relative
generation $dT+L$, to see whether $w$ has property $\ca_T$. If so, put
$w$ into $S_1$ and the vertices in $D(w)$ into $S_2$.  If $S_1$ is
large enough to guarantee that $v$ has $\ca$, stop. (At this stage,
this happens if and only if $|S_1|\ge k$.) Otherwise, examine each $w$
in generation 2 that is not a descendant of a vertex in $S_1$ in
turn. For each $w$, reveal its descendants up to $dT+L$ generations later,
and test whether $w$ has $\ca_T$. If so, include $w$ into $S_1$ and
the vertices of $D(w)$ into $S_2$.  Continue in this way until either
(a) the set $S_1$ shows that $\ca$ holds (i.e.,
$S_1$ contains a set of targets witnessing $\brs\of\ca_T$), or (b) we have
tested all vertices in generations up to $2L+1$ not descended from
vertices in $S_1$ for property $\ca_T$ without finding such a set. In
case (b), $\ca=\brs\of \ca_T$ does not hold. Hence, $v$ has $\ca$ if
and only if (a) holds.

Let us condition on $v$ having property $\ca$.  The key points
are that this event is guaranteed by the vertices in $S_2$ and their
ancestors, and that we have not examined the children of any vertex in
$S_2$. Let us suppose that at most $n^{1/10}$ vertices have been examined,
an event of probability $1-o(n^{-1})$ by \eqref{small}. We now examine
the vertices $w\in S_2$ one by one. For each $w$ we explore its descendants
for the next $1+dT+L$ generations to test whether, in the branching process
rooted at $v$, the particle $w$ has the property $\br_1\of\ca_T$. We abandon
the exploration associated to a given $w$ if we reach more than $n^{1/10}$
vertices in this exploration. If the test is successful, we mark $w$.
Note that, from \eqref{small}, the probability
that any exploration is abandoned is $o(n^{-1})$. 
Provided no previous exploration has been abandoned, the argument above
shows that we may couple the descendants of $w$ to agree with a
branching process $X_{\la'}$ with probability $1-o(1)$. Hence,
if $n$ is large enough,  the probability
that we mark $w$ is at least
\[
 \Pr_{\la'}(\br_1\of\ca_T)-o(1) \ge \beta(\la_2),
\]
using \eqref{a2}.
In summary, ignoring an error probability of $1-o(n^{-1})$, we can view each
$w\in S_2$ as marked independently with probability (at least) $\beta(\la_2)$.

Let $\cT_v$ be the tree rooted at $v$ with leaves the marked vertices in $S_2$.
Since $T/\log\log n\to \infty$, we have $k^{-(2^T+2)d}\le n^{-10}$ if 
$n$ is large enough. Hence,
from \eqref{c1}, with probability $1-o(n^{-1})$ every vertex of $S_1$ has property $\ev_T$
guaranteed by its descendants in $\cT_v$, so with probability $1-o(n^{-1})$
the tree $\cT_v$ has the property $\cB=\brs\of\ev_T$ when viewed as a branching process. Let
us suppose that $\cT_v$ has property $\cB$, and let $\cT_v'$ be a minimal subtree with this
property. Recalling that $\ev_T=\br_{Td}\of\rb_L$,
note that every non-leaf
of $\cT_v'$ has degree at least $k$ in $\cT_v'$.

We claim that each leaf $w$ of $\cT_v'$ has property $\ca$ in the graph $\tG$. To see this,
let $x$ be the ancestor of $w$ in $\cT_v'$ that is $L$ generations above $w$,
so the subtree $\cT'$ of $\cT_v'$ rooted at $x$ is a tree of height
$L$ that is minimal with respect to having property $\rb_L$.
All leaves of $\cT'$ are leaves of $\cT_v'$ and hence marked, so they have
the property $\br_1\of\ca_T$ in the tree $G_v[2s]$ rooted at $v$.
Hence we may take $k-1$ children of each leaf and add them to $\cT'$ to form a tree $\cT$
with property $\rb_L\of\br_1$ in which every leaf has property $\ca_T$.
By Lemma~\ref{brp}, taking $w$ as the root, this tree $\cT$ has a subtree
$\cW$ with property $\Ms$ all of whose leaves are leaves of $\cT$
and hence have property $\ca_T$. Thus $w$ has property $\brs\of\ca_T=\ca$ in the graph,
as claimed.

We have shown that with probability $1-no(n^{-1})=1-o(1)$ {\em every}
vertex $v$ with property $\ca$ is the root of a tree $\cT_v'$ in which each non-leaf
has degree at least $k$ and each leaf has property $\ca$ (in the graph $\tG$).
When this happens, the union of the trees $\cT_v'$ is a subgraph of $\tG$ with minimum degree
at least $k$ containing all $N$ vertices with property $\ca$, so
$c_k(G(n,\la/n))\ge N$. This completes our proof of the
lower bound.

\medskip
The proof of the upper bound on $c_k(G(n,\la/n))$ is much simpler:
the events $\brp_d$ decrease to $\brp$ as $d\to\infty$.
Hence, given any $\eps>0$, there is a $d$ with $\Prl(\brp_d)\le \bp(\la)+\eps$.
Exploring the first $d$ neighbourhoods of a random vertex $v$ of $G(n,\la/n)$ as above,
the probability that we encounter a cycle is $o(1)$. Let us say that $v$ has property
$\brp_d$ if its $d$-neighbourhood is a tree with property $\brp_d$.
Considering the local exploration described above for $\tG$,
in the simpler context of $G(n,\la/n)$ it is easy to check that
$\Prl(\brp_d)n+o_p(n)$ vertices $v$ have property $\brp_d$, while $o_p(n)$ vertices
 have a cycle in their $d$-neighbourhood. But any vertex in the $k$-core must have one of these
properties, so
\[
 c_k(G(n,\la/n))\le \Prl(\brp_d)n+o_p(n) \le \bp(\la)+\eps n+o_p(n).
\]
As $\eps>0$ was arbitrary, this completes the proof.
\end{proof}

Theorem \ref{t1} implies that the natural coupling between the neighbourhoods
of a vertex $v$ in $G(n,\la/n)$ and the branching process $X_\la$ can be adapted
to the $k$-core. With $k\ge 2$ fixed as usual, let $C_k(G)$ denote the $k$-core
of $G$, and let $X_\la'$ denote the (possibly empty) branching process
consisting of all particles in $X_\la$ that are in a $k$-regular tree in $X_\la$
containing the initial particle.

\begin{corollary}\label{cor1}
Let $k\ge 2$, $\la>0$ and $L$ be fixed. If $\beta(\la)$ is continuous at $\la$,
then the first $L$ neighbourhoods
in $C_k(G(n,\la/n))$ of a random vertex of $G(n,\la/n)$ may be coupled
with the first $L$ generations of $X_\la'$ so as to agree \whp.
\end{corollary}
\begin{proof}
We use only the result of Theorem \ref{t1}, not anything from the proof.
We may assume that $\beta(\la)>0$, as otherwise $C_k=C_k(G(n,\la/n))$ is empty \whp,
while $X_\la'$ is empty with probability $1$.
Since $\brp_d\downto \brp$, if $d(n)\to\infty$ then $\Pr_\la(\brp_{d(n)})=\Pr_\la(\brp)+o(1)$.
Note that $\brp_{d(n)}\supset \brp$.

If $d(n)\to\infty$ sufficiently slowly then, by standard arguments, we may couple
the branching process $X_\la$ with the neighbourhood process $Y=Y(v,G(n,\la/n))$
so that \whp\ they agree for the first $d(n)$ generations. (This follows from the
$d$ constant case of the coupling.)
If $v$ is in the $k$-core $C_k$, then either $\brp_{d(n)}$ holds for $Y$,
or the first $d(n)$ neighbourhoods of $v$ contain a cycle.
Now
\[
 \Pr(Y\in \brp_{d(n)}) = \Pr(X_\la\in \brp_{d(n)}) +o(1) = \Pr(X_\la\in \brp)+o(1)
 = \Pr(v\in C_k) +o(1),
\]
where the first step is from the coupling and the last
from Theorem \ref{t1}.
Provided $d(n)=o(\log n)$, which we may assume, $o_p(n)$ vertices
$v$ have a cycle in their $d(n)$-neighbourhoods. Thus, \whp\ either all three
of $\{v\in C_k\}$, $\{Y\in \brp_{d(n)}\}$ and $\{X_\la\in \brp\}$ hold, or none.

The argument above holds with $d(n)-L$ in place of $d(n)$. Let us call a vertex
$w$ {\em exceptional} if it is within distance $d(n)-L$ of a cycle, 
or if its neighbourhoods to distance $d(n)-L$ have the property
corresponding to $\brp_{d(n)-L}$ but $w$ is not in the $k$-core. 
All other vertices are {\em normal}. Then the probability
that a given vertex $w$ is exceptional is $o(1)$, so there are $o_p(n)$ exceptional
$w$. Since, for any fixed $t$, the $t$-neighbourhood of any set of $o_p(n)$ vertices
in $G(n,\la/n)$ has size $o_p(n)$, returning to the random vertex $v$,
with probability $1-o(1)$ all vertices within distance $L$ of $v$ are normal.

Suppose that the coupling above succeeds for the neighbourhoods of a random vertex $v$
to distance $d(n)$, and let $S$ be the set of vertices near (i.e., within distance $L$ of) $v$
that are in the $k$-core. As \whp\ every vertex near $v$ is normal, \whp\ $S$ is
the set of vertices corresponding to particles $y$ in the first $L$ generations of $X_\la$
with the property $\brp_{d(n)-L}$ in the tree $X_\la$ with $y$ as root.
As $L$ is fixed and $d(n)\to \infty$, this latter set agrees \whp\ with the set $Z$ of
particles  $z$
in the first $L$ generations of $X_\la$ having property $\brp$ in $X_\la$
with $z$ as root, i.e., contained
in some infinite $k$-regular tree. But the set of particles in $Z$ reachable from the root
$x_0$ by a path in $Z$ is just the first $L$ generations of $X_\la'$, so \whp\ the
first $L$ generations of $X_\la'$ correspond to the first $L$ neighbourhoods
of $v$ in the $k$-core, as required.
\end{proof}

The argument above shows that,
once we have the asymptotic number of vertices in the $k$-core $C_k$, we can count
up to $o_p(n)$ all `local' structures in $C_k$, including, for example, the
number of vertices of a given degree. Of course, we can also estimate the sum
of the degrees, say, and hence the number of edges. More generally,
we can estimate the sum of any function of the $L$-neighbourhood with a well-behaved tail,
where the contribution from the $o(1)$-probability case when the coupling breaks
down is small.

Let $f(v,G)$ be an isomorphism-invariant function of graphs $G$ rooted at a vertex $v$.
As in \cite{kernels}, we call $f$ an $L$-{\em neighbourhood function} if it depends
only on the subgraph of $G$ induced by the vertices within distance $L$ of $v$.
(Thus the degree of $v$ in $G$ is a $1$-neighbourhood function.)
We interpret $f(X_\la)$ in the natural way, viewing the branching process as a tree
with root $x_0$.

\begin{theorem}\label{tcount}
Let $k\ge 2$, $\la\ne\lac(k)$, $L\ge 1$ and an $L$-neighbourhood function $f$ be given. 
If $f(v,G)$ is bounded by a polynomial of the number of vertices in the $L$-neighbourhood
of $v$, then
\[
 S_n\=\frac{1}{n}\sum_{v\in C_k} f(v,C_k) \pto \E f(X_\la'),
\]
where $C_k$ denotes the $k$-core of $G(n,\la/n)$.
\end{theorem}

The simple proof follows exactly that of the corresponding result
for the $2$-core of a more general graph, Lemma 11.11 of \cite{kernels}.
As the result is very unsurprising in the light of Corollary \ref{cor1}, we omit the proof.
The condition that $f$ be polynomially bounded can be replaced by a fourth moment
condition, as in \cite{kernels}.

\section{Inhomogeneous random graphs}\label{sec_nonunif}

Many random graph models have been considered in which edges are independent but different
possible edges have different probabilities,
including numerous `sparse' models generalizing $G(n,\la/n)$, where
the expected number of edges is linear in the number of vertices; see
\cite{kernels} and the references therein.
A very general model of this type was introduced by Bollob\'as, Janson and Riordan
in \cite{kernels}; let us recall the definitions.

A {\em ground space} is a pair $(\sss,\mu)$, where
$\sss$ is a separable metric space and $\mu$ is a Borel probability
measure on $\sss$. Mostly, we shall consider the cases $\sss$ finite,
with $\mu$ any (strictly positive) measure on $\sss$, and
$\sss=(0,1]$ or $\sss=[0,1]$, with $\mu$ Lebesgue measure.

A set $A\subseteq\sss$ is a \emph{\mucs}
if $A$ is (Borel) measurable and $\mu(\ddd A)=0$, where
$\ddd A$ is the boundary of $A$.

A {\em vertex space} $\vxs$ is a triple $(\sss,\mu,\xss)$, where
$(\sss,\mu)$ is a ground space and, for each $n\ge 1$, $\xs$
is a random sequence $(\xn_1,\xn_2,\ldots,\xn_n)$ of $n$ points of $\sss$,
such that
\begin{equation}
  \label{a2a}
\nu_n(A)\=\#\{i:\xn_i\in A\}/n \pto \mu(A)
\end{equation}
as $n\to\infty$, for every \mucs\ $A$.

The sequence $\xn_i$ will give the {\em type} of each vertex $i$ in a graph on $[n]$
still to be defined. The convergence condition \eqref{a2a} says that distribution
of the types of the vertices is essentially $\mu$. More precisely, the empirical
distribution $\nu_n$ of the types converges in probability to $\mu$. In the {\em finite-type}
case, where $\sss=[r]=\{1,2,\ldots,r\}$ and $\mu\{i\}>0$ for $i=1,2,\ldots,r$,
the condition says that the number $n_i$ of vertices of each type $i$
satisfies $n_i/n\pto \mu\{i\}$. No assumption is made about the dependence between
the types of different vertices.

A \emph{kernel} $\kappa$ on a ground space $(\sss,\mu)$
is a symmetric non-negative (Borel) measurable
function on $\sss\times \sss$.
A kernel on a vertex space $(\sss,\mu,\xss)$ is simply a kernel on
$(\sss,\mu)$.

Finally, let $\kappa$ be a kernel on the vertex space $\vxs$. Given
the (random) sequence $(\xn_1,\dots,\xn_n)$, let $\gnkx$ be the random
graph on $[n]$ in which edges are present independently,
and the probability that the edge $ij$ is present is
\begin{equation}
  \label{pij}
 \pij\=\min\big\{\kk(\xn_i,\xn_j)/n,1\big\}.
\end{equation}

In other words, we assign vertices types from $\sss$ with the types
asymptotically distributed according to $\mu$, and then join
vertices of types $x$ and $y$ with probability $\kk(x,y)/n$.
If $\kk$ takes the constant value $\la$, then $\gnkx$ is exactly
$G(n,\la/n)$.

The model just described includes many models of inhomogeneous
graphs previously defined. As the model is very general, it may
help to bear in mind a few special cases. The most fundamental
is the finite-type case described above. This is an extremely natural
generalization of $G(n,\la/n)$, and has been considered
earlier by several authors (e.g., S\"oderberg \cite{Sod1}),
perhaps with slightly different assumptions
on the distribution of the types. (For example, the types
of the vertices being independent, or $n_i$ equal to $\lfloor\mui n\rfloor$ 
or $\lceil \mui n\rceil$.)
Several other interesting special cases
have $\sss=(0,1]$, $\mu$ Lebesgue measure, and $\xn_i=i/n$,
so the probability that the edge $ij$ is present is $\pij=\kk(i/n,j/n)/n$.
Taking $\kk(x,y)=1/\max\{x,y\}$ gives $\pij=1/\max\{i,j\}$;
minor variants of the corresponding random graph were
introduced independently by Dubins in 1984 (see \cite{KW})
and by
Callaway, Hopcroft, Kleinberg, Newman and Strogatz~\cite{CHKNS} in 2001;
see \cite{kernels}.

In order to prove results about $\gnkx$, some additional assumptions are needed
to avoid pathologies. Following \cite{kernels}, we assume throughout that
$\kk$ is continuous almost everywhere on $\sss\times\sss$, that
\begin{equation}\label{int}
 \int_{\sss\times\sss} \kk(x,y)\dd\mu(x)\dd\mu(y)<\infty,
\end{equation}
and that
\begin{equation}\label{nopath}
 \frac1n  \E e\bb{\gnkx} \to \frac12\iint_{\sss^2}\kk(x,y)\dd\mu(x)\dd\mu(y)
\end{equation}
as $n\to \infty$.
The last condition says that the number of edges of $\gnkx$ is `what it should be',
at least in expectation.
Without this condition, it is impossible to relate the behaviour of $\gnkx$ to
that of $\kk$ (or rather, to that of the branching process $X_\kk$ defined below):
changing $\kk$ on a set of measure zero should not affect the model, but $\gnkx$
may depend {\em only} on the values of $\kk$ on such as set, for example, when
$\xn_i=i/n$.
Surprisingly, the very natural and rather weak
condition \eqref{nopath} is enough to enable many results about $\gnkx$ to be proved;
see \cite{kernels} for a discussion of this.

We shall need one further definition: a kernel $\kk$ on a ground space $(\sss,\mu)$
is \emph{reducible} if
\[
  \text{$\exists A\subset \sss$ with
$0<\mu(A)<1$ such that $\kappa=0$ a.e.\ on $A\times(\sss\setminus
A)$},
\]
and {\em irreducible} otherwise. Roughly speaking, reducible kernels
correspond to disconnected graphs; much of the time, nothing is lost
by considering only irreducible kernels.

The actual definitions in \cite{kernels} are slightly more general in two
ways; firstly, the kernel $\kk$ is allowed to depend somewhat on $n$. We shall
not bother with this additional generality here. Secondly, the number
of vertices of $\gnkx$ need not be exactly $n$, but may be $n+o_p(n)$.
This latter relaxation complicates only the notation, not the proofs.

The key to the analysis of $\gnkx$ turns out to be the multi-type
Galton--Watson branching process $X_\kk$ associated to $(\sss,\mu,\kk)$.
This starts with a single particle $x_0$ whose type is distributed
according to $\mu$. Each particle of type $x$ has a set of children
whose types are distributed as a Poisson process on $\sss$ with
intensity $\kk(x,y)\dd\mu(y)$. In other words, the number of children with
types in a subset $A\subseteq\sss$ has a Poisson distribution with
mean $\int_A \kk(x,y)\,d\mu(y)$, and these numbers are independent
for disjoint sets $A$ and for different particles; see, for example,
Kallenberg \cite{Kall}.
Sometimes it will be convenient to start the process with a particle
of a fixed type $x$ instead of a random type. We write
$X_\kk(x)$ for this branching process.

If $\kk$ takes the constant value $\la$, then $X_\kk$ is just $X_\la$ as
defined in Section~\ref{sec_unif}. In general, the relationship of $X_\kk$ to $\gnkx$
is the same as that of $X_\la$ to $G(n,\la/n)$.
For example, under suitable regularity conditions, which certainly
hold in the finite-type case, the first few neighbourhoods of a random
vertex $v$ of $\gnkx$ may be coupled with the branching process
$X_\ka$ in the natural sense (so that the type of a vertex is the same
as the type of the corresponding particle in the branching process); see \cite{kernels}.
One of the main results of \cite{kernels} is that, under suitable mild
assumptions, the size of the giant component in $\gnkx$ is asymptotically
$n$ times the probability that the branching process $X_\ka$ never dies
out. This generalizes the classical result giving the size of the giant
component of $G(n,\la/n)$ in terms of $X_\la$.

We shall write $\Prk$ for the probability measure associated to $X_\kk$,
and $\Prkx$ for that associated to $X_\kk(x)$. Thus,
if $\ca$ is some property of rooted trees, we write
$\Prk(X_\kk\in \ca)$, or simply $\Prk(X\in \ca)$, for the probability
that $X_\kk$ has this property when viewed as a tree.

Turning to the $k$-core, let us define $\beta(\kk)$ and $\bp(\kk)$
as before, but for $X_\kk$. Thus $\beta(\kk)=\Prk(\br)$,
and $\bp(\kk)=\Prk(\brp)$. We shall
also write $\beta_x(\kk)$ for $\Prkx(\br)$, and $\bp_x(\kk)$
for $\Prkx(\brp)$. By analogy with the result for
$G(n,\la/n)$, we expect the $k$-core
of $\gnkx$ to have size $\bp(\kk)n+o_p(n)$, at least
under suitable conditions.

\subsection{The finite-type case}

The strategy of starting with the finite-type case, and then
using approximation and monotonicity arguments to attack the general case,
is used throughout \cite{kernels}. We use the same strategy here; 
we have written the proof of Theorem \ref{t1} so that
it carries over almost immediately to the finite-type case.
Recall that we write $c_k(G)$ for $|C_k(G)|$, where $C_k(G)$ is
the $k$-core of a graph $G$.

\begin{theorem}\label{t2}
Let $k\ge 2$ be fixed.
Let $\ka$ be a kernel on a vertex space $\vxs=(\sss,\mu,\xss)$ with
$\sss=[r]$, $r\ge 1$, and $\mui>0$ for each $i$.
If the function $\la\mapsto\beta(\la\ka)$ is continuous
at $\la=1$, then
\[
 c_k(\gnkx) = \bp(\ka)n+o_p(n).
\]
\end{theorem}

\begin{proof}
We may assume without loss of generality that $\kk$ is irreducible; indeed,
in the finite-type case, if $\kk$ is reducible the graph $\gnkx$ may be written
as the disjoint union of two or more graphs given by instances
of the same model with irreducible kernels. (As each $n_i$ is random,
these graphs have a random number
of vertices, but this does not matter.)

As in the proof of Theorem \ref{t1}, the upper bound is easy. Indeed,
arguing as for $G(n,\la/n)$, it
is easy to see that for fixed $d$ one can couple the $d$-neighbourhoods of a random
vertex $v$ of $\gnkx$ with the first $d$ generations of $X_\ka$ to agree with
probability $1-o(1)$. The upper bound on the $k$-core follows as in the uniform case.

For the lower bound, it again suffices to prove that if $\la<\la'$, then
\whp\ the $k$-core of $\gnxx{\la'\ka}$ contains at least $\bp(\la\ka)$ vertices.
The proof is essentially the same as in the uniform case, {\em mutatis mutandis};
we indicate the changes briefly.

The proof of Lemma \ref{l1} extends unchanged. 
Since $\mu\{i\}>0$ for each $i$,
it then follows that
\[
 \Pr_{\la\ka,i}(\rb_d\of\br)\upto \beta_i(\la\ka)
\]
as $d\to\infty$, for each $i\in \sss=[r]$. For $\vp=(p_1,\ldots,p_r)$, 
let $r_i(\la\ka,d,\vp)$ be the probability that $\rb_d\of M$
holds if we mark the particles in generation $d$ of $X_{\la\ka}(i)$ independently,
marking a particle of type $j$ with probability $p_j$. In place of Lemma \ref{l2} we obtain
\[
 r_i(\la_2\ka,d,\vp) \ge \beta_i(\la_1\ka)
\]
for every $i$, where $p_j=\beta_j(\la_1\ka)$; in the proof,
we add $\Po\bb{(\la_2-\la_1)\ka(i,j)}$ `extra' particles of some type $j$ with $\ka(i,j)>0$.
The irreducibility of $\kk$ guarantees that
such a type $j$ exists, and also that $\beta_j(\la_1\kk)>0$ for every
$j$ whenever $\beta(\la_1\kk)>0$.

Lemma \ref{l_small} adapts immediately, as Harris' Lemma applies to random subsets
obtained by selecting elements independently, even if the selection probabilities are
different for different elements.
Next, in place of \eqref{close} we have
\[
 \pgem\left(\la'\sum_j \ka(i,j)\beta_j(\la_1\ka)\right) \ge \beta_i(\la_2\ka)
\]
for each $i$. Arguing as for \eqref{pl}, \eqref{cpl} but taking
\[
 \eta=\min\bigl\{k^{-3d},\,
  \min_i \bigl\{ \bb{\Pr_{\la_2\ka,i}(\rb_{d_1}\of\br)-\beta_i(\la_1)}\mui\,\bigr\} \bigr\} >0,
\]
we find an event $\el$ depending on the first $L$ generations with
\[
 \Pr_{\la_2\ka}(\rb_L\of\br\mid X[L]) \ge 1-k^{-3d}
\]
whenever $\el$ holds, and $\Pr_{\la_2\ka,i}(\el)\ge \beta_i(\la_1\ka)$ for each $i$.
It follows as before that $\Pr_{\la_2\ka,i}(\ca_t)\ge \beta_i(\la_1\ka)$ for each $i$,
and the rest of the proof is essentially unchanged.
\end{proof}

\begin{remark}
The continuity condition in Theorem \ref{t2} is non-trivial: unlike in
the uniform case, $\la\mapsto \beta(\la\ka)$ may have several discontinuities.
One way this can happen is when $\kk$ is reducible,
so the graph decomposes into two separate pieces, whose
$k$-cores may emerge at different points. Perhaps surprisingly, this is not the only way.
Let $k=3$, and set $\sss=\{1,2\}$ with $\mu\{1\}=\mu\{2\}=1/2$.
Let $\ka(1,1)=2000$, $\ka(2,2)=2$ and $\ka(1,2)=\ka(2,1)=1/100$, say.
It is easy to check that a $k$-core first emerges near $\la=\lac/1000$,
where $\lac$ is the critical parameter for the emergence of a $3$-core in $G(n,\la/n)$;
at this point, the vertices of type $1$ form a uniform random graph with large enough
average degree to contain a $k$-core.
When $\la$ is close to $\lac$, the probability $p_2=\beta_2(\la\kk)$ that a vertex
of type $2$ is in the $k$-core is related to the largest solution
to $p=\pgem(\la p+\la p_1/200)$, where $p_1$ is the (unknown) probability that
a vertex of type $1$ is in the $k$-core. Since $0<p_1<1$, it is easy to check
that the solution jumps near $\la=\lac$, and in fact that
$p_2$ jumps from around $\pgek(\lac/100)$ to around $\beta(\lac)$.
\end{remark}

As in the uniform case, Theorem \ref{t2} implies the equivalents
of Corollary \ref{cor1} and Theorem \ref{tcount} for the finite-type
case. The statements and proofs are direct translations of those in 
Section \ref{sec_unif}, so we omit them.

\subsection{The general case}

Theorem \ref{t2} extends easily to more general kernels under some mild assumptions.
To state these, we need another definition:
a kernel $\kk'$ is {\em regular finitary} if there is a partition of $\sss$
into a finite number $r$ of \mucs s $S_i$, such that $\kk'$ is constant
on each $S_i\times S_j$. Such kernels correspond to finite-type kernels in an obvious way.

When studying the $k$-core in $G(n,\la/n)$, we assumed that $\beta(\cdot)$
(or, equivalently, $\bp(\cdot)$) was continuous from below at $\la$, i.e.,
that $\beta(\la_n)\upto \beta(\la)$ for every sequence $\la_n\upto \la$.
Of course, it makes no difference if we consider only one sequence $\la_n$,
as long as this increases strictly to $\la$. It is this latter, formally weaker,
condition that we shall adapt to general kernels, with the restriction
that our kernels tending up from below be of finite type.

When studying $\gnkx$,
we shall assume that the functional $\kk\mapsto\bp(\kk)$ is {\em continuous from below at $\kk$},
in the weak sense that
\begin{equation}\label{cts}
 \bp(c_m\kk_m)\to\bp(\kk)
\end{equation}
for {\em some} sequence $\kk_m$ of regular finitary kernels with $\kk_m(x,y)\le \kk(x,y)$
for all $x$ and $y$, and some sequence $c_m$ of real numbers with $c_m<1$ for all $m$.
(Equivalently, we require $\bp(\kk_m)\to\bp(\kk)$ for some finite-type kernels $\kk_m$
with $\sup_{x,y}\kk_m(x,y)/\kk(x,y) <1$ for all $m$.
When $|\sss|=1$, this is equivalent to continuity from below in the usual sense.)
As we have seen, this is
a non-trivial condition; however, we expect it to hold for almost all $\kk$,
in some imprecise sense. We return to this later. Also, it may well be the
case that, if $\kk$ is irreducible, then
this condition is equivalent to the condition that $\bp(\kk_m)\to\bp(\kk)$
for {\em all} sequences $\kk_m$ of (arbitrary) kernels increasing to $\kk$ (almost everywhere).
Note that $\bp$ is always continuous from above:
if $\kk_m$ is a sequence of kernels decreasing pointwise to $\kk$,
then one can couple the branching processes $X_{\kk_m}$ and $X_\kk$ so that
$X_{\kk_1}\supset X_{\kk_2}\supset \cdots$, and $X_\kk=\bigcap_m X_{\kk_m}$.
It follows that $\bp(\kk_m)\downto \bp(\kk)$.

Our continuity assumption can be viewed as the assertion that $o(n)$ edges will
not change the size of the $k$-core much. Due to the flexibility of the model,
it would be unreasonable to expect a precise result without this condition;
as shown in \cite{kernels}, the condition \eqref{nopath} used to exclude
pathologies still permits the insertion of $o_p(n)$ edges into $\gnkx$
in a more or less arbitrary way, while changing $\kk$ on a set of measure
zero, which does not alter $\bp(\kk)$.

Our second assumption will state essentially that a small set of `exceptional'
vertices (or edges) cannot have many vertices in its $L$-neighbourhood, for any fixed $L$.
Perhaps the most natural form of this assumption involves counting paths:
let $P_\ell(G)$ denote the number of $\ell$-edge paths in a graph $G$, let
\[
 \alpha_\ell(\kk) \=
\frac{1}{2}
 \int_{\sss^{\ell+1}} \kk(x_0,x_1) \kk(x_1,x_2) \dotsm
 \kk(x_{\ell-1},x_\ell)\dd \mu(x_0)\dotsm \dd \mu(x_\ell),
\]
and suppose that
\begin{equation}\label{pconv}
 \frac{1}{n}\E\left(P_\ell\bb{\gnkx}\right) \to \alpha_\ell(\kk) <\infty
\end{equation}
as $n\to\infty$, for each $\ell\ge 1$. The convergence condition says essentially
that the expected number of paths in $\gnkx$ is `what it should be'.
As shown in \cite{kernels}, it holds whenever $\kk$ is bounded, or whenever
the types of the vertices are independent with distribution $\mu$.
Theorem 17.1 in \cite{kernels} shows that convergence in
expectation in \eqref{pconv}
implies convergence in probability; we shall not directly use this fact.

By Lemma 7.3(ii) of \cite{kernels}, there is a sequence of regular finitary
kernels $\kk_m$ with $\kk_m\le \kk$ pointwise, such that $\kk_m(x,y)\upto\kk(x,y)$ for almost
all $(x,y)\in \sss\times\sss$.
From dominated convergence we have $\alpha_\ell(\kk_m)\upto \alpha_\ell(\kk)$ for each $\ell$.
It is easy to check that for a regular finitary kernel $\kk'$ we have
\[
 \frac{1}{n}\E\left(P_\ell\bb{\gnxx{\kk'}}\right) \to \alpha_\ell(\kk')
 \hbox{ and }
 \frac{1}{n}P_\ell\bb{\gnxx{\kk'}} \pto \alpha_\ell(\kk').
\]
In particular, if \eqref{pconv} holds then as $m\to\infty$ we have
\begin{equation}\label{asum2}
 \lim_{n\to\infty}\frac{1}{n}\E\left(P_\ell\bb{\gnxx{\kk_m}}\right) \to  \lim_{n\to\infty}\frac{1}{n}\E\left(P_\ell\bb{\gnxx{\kk}}\right) <\infty.
\end{equation}
Since $\kk_m\le \kk$, we may couple $\gnxx{\kk_m}$ and $\gnkx$ so that $\gnxx{\kk_m}\subset \gnkx$.
Condition \eqref{asum2} says essentially that, if $\ell$ is fixed and $m(n)\to\infty$, then
almost all paths of length $\ell$ in $\gnkx$ are already present in the subgraph $\gnxx{\kk_{m(n)}}$.

The actual assumption we shall need is a tiny bit weaker. We shall assume that there is
an increasing sequence $\kk_m$ of regular finitary kernels with  $\kk_m\le \kk$ pointwise,
and a coupling $\gnxx{\kk_m}\subset \gnkx$, with the following property:
for any $\ell\ge 1$ and any $\eps>0$ there is an $m_0=m_0(\ell,\eps)$ such that
\begin{equation}\label{asum3}
 \Pr\left( |V_{\ell,m}|\ge \eps n\right) \le \eps
\end{equation}
for all $m\ge m_0$ and all large enough $n$,
where $V_{\ell,m}$ is the set of vertices that are endpoints of a path
of length $\ell$ in $\gnkx$ not present in $\gnxx{\kk_m}$. In other words,
if $m(n)\to\infty$, then $|V_{\ell,m(n)}|=o_p(n)$. As a path has a bounded number of endpoints,
\eqref{asum2} immediately implies this condition: all paths present
in $\gnxx{\kk_m}$ are also present in $\gnkx$, and the expected number of additional paths is small.

Assuming continuity and \eqref{asum3}, it is very easy to extend Theorem \ref{t2}.

\begin{theorem}\label{t3}
Let $k\ge 2$ be fixed.
Let $\kk$ be a kernel on a vertex space $\vxs$. Suppose that $\kk$ is continuous
almost everywhere on $\sss\times\sss$, and that \eqref{int} and \eqref{nopath} hold.
If, in addition, $\kk\mapsto\bp(\kk)$ is continuous at $\kk$ in the sense of \eqref{cts}, and
\eqref{asum3} holds, then
\[
 \frac{1}{n} c_k\bb{\gnkx} \pto \bp(\kk)
\]
as $n\to\infty$.
\end{theorem}

\begin{proof}
Let $\kk_m$ be a sequence of regular finitary kernels with $\kk_m\le \kk$ pointwise,
and $c_m$ a sequence of real numbers with $c_m<1$, 
such that $\bp(c_m\kk_m)\to \bp(\kk)$. Such sequences exists by our
continuity assumption \eqref{cts}. Let $\eps>0$ be arbitrary.

Since $\bp(c_m\kk_m)\to \bp(\kk)$, there
is an $m$ with $\bp(c_m\kk_m)\ge \bp(\kk)-\eps$.
Since $c\mapsto \bp(c\kk_m)$ is an increasing, continuous function of $c$, 
there is a $c$ with $c_m<c<1$ at which this function is continuous.
As $c\kk_m(x,y)\le \kk_m(x,y)\le \kk(x,y)$ for all $x$ and $y$, we may couple $\gnxx{c\kk_m}$
and $\gnkx$ so that $\gnxx{c\kk_m}\subset \gnkx$ for every $n$.
Applying Theorem \ref{t2} to the finite-type kernel corresponding to $c\kk_m$
in the natural way, we have
\[
 c_k\bb{\gnkx} \ge c_k\bb{\gnxx{c\kk_m}} = \bp(c\kk_m)n+o_p(n),
\]
so \whp
\[
 c_k\bb{\gnkx}\ge (\bp(c\kk_m)-\eps)n \ge (\bp(c_m\kk_m)-\eps)n \ge (\bp(\kk)-2\eps)n.
\]

To prove the upper bound, let $\kk_m\le \kk$ be a sequence
of regular finitary kernels satisfying \eqref{asum3}. (We shall
not assume that $\bp(\kk_m)\to \bp(\kk)$ for this sequence.)
Since $\brp_d\downto \brp$, there is a $d$ with
$\Pr_\kk(\brp_d)\le \bp(\kk)+\eps$.
By \eqref{asum3},
for any $\eta>0$ there is an $m$ such that
\begin{equation}\label{vsmall}
 \Pr\left(|V_{\ell,m}|\ge \eps n/d \right) \le \eta/d
\end{equation}
for all $1\le \ell\le d$ and all large enough $n$.

Let $v$ be a vertex of $\gnkx$, so $v$ is also a vertex of $\gnxx{\kk_m}$.
If $v$ is in the $k$-core of $\gnkx$, then either the $d$-neighbourhood $\Ga_d(v,\gnkx)$ of $v$
in $\gnkx$ contains a tree with the property $\brp_d$,
or it contains a cycle.
Thus, one of the following three cases holds: (i) $\Ga_d(v,\gnxx{\kk_m})$ contains a tree with property $\brp_d$,
(ii) $\Ga_d(v,\gnxx{\kk_m})$ contains a cycle, or (iii)
$\Ga_d(v,\gnkx)$ contains an edge not present in $\gnxx{\kk_m}$.
In case (iii), $v$ is an endpoint of a path of length at most $d$
contained in $\gnkx$ but not in $\gnxx{\kk_m}$, i.e., $v\in \bigcup_{\ell\le d}V_{\ell,m}$.
From \eqref{vsmall}, with probability $1-\eta$ there are at most $\eps n$ such $v$.
By the finite-type equivalent of Corollary \ref{cor1}, the number of vertices
with property (i) is
\[
 \Pr_{\kk_m}(\brp_d)n +o_p(n) \le \Pr_{\kk}(\brp_d)n +o_p(n) \le \bp(\kk)n+\eps n +o_p(n).
\]
Also, by the same result, or by directly counting short cycles (noting that $\kk_m$ is bounded),
only $o_p(n)$ vertices $v$ have property (ii).
It follows that, for any $\eta$, if $n$ is large enough we have
$c_k(\gnkx)\le (\bp(\kk)+3\eps)n$ with probability at least $1-\eta-o(1)$.

As $\eta>0$ was arbitrary, we have shown that \whp{}
\[
 (\bp(\kk)-2\eps)n \le c_k(\gnkx) \le (\bp(\kk)+3\eps)n.
\]
As $\eps>0$ was arbitrary, this completes the proof.
\end{proof}

As noted above, the rather cumbersome assumption \eqref{asum3} holds whenever
\eqref{pconv} holds. In the next section, we shall consider an interesting example
which satisfies \eqref{asum3} but not \eqref{pconv}.

In some sense, under mild assumptions, Theorem \ref{t3} is the complete
answer to the question `how large is the $k$-core in $\gnkx$?' In another
sense, it is just the beginning: it remains to determine $\bp(\kk)$.
Although there is no further combinatorics involved, depending on the form
of the solution required, this is in general very difficult. However,
it is very easy to determine $\bp(\kk)$ in terms of the solution to a
certain functional equation.

{}From the definition of $\beta_x$ and the definition 
of the branching process, we have
\begin{equation}\label{fe}
 \beta_x(\kk) = \pgem\left(\int_\sss \kk(x,y)\beta_y(\kk) \dd\mu(y)\right).
\end{equation}
Indeed, starting with a particle of type
$x$, the number of children having property $\br$
is Poisson with mean $\int \kk(x,y)\beta_y(\kk)\dd\mu(y)$.
With $\kk$ fixed, we regard \eqref{fe} as a functional equation in a function
$x\mapsto \beta_x(\kk)$. Arguing as in the uniform case, it is easy
to see that $\beta_x(\kk)$ is given by the maximum solution to this equation
(i.e., the supremum of all solutions, which is itself a solution).
Arguing as for \eqref{fe}, we have
\begin{equation}\label{bp1}
 \bp_x(\kk) = \pgek\left(\int_\sss \kk(x,y)\beta_y(\kk) \dd\mu(y)\right),
\end{equation}
while, of course,
\begin{equation}\label{bp2}
 \bp(\kk) = \int_\sss \bp_x(\kk) \dd \mu(x).
\end{equation}

Turning to the continuity of $\kk\mapsto \bp(\kk)$ at $\kk$,
suppose that $\kk_m\upto \kk$, and define $\gamma_x$ by
\begin{equation}\label{gdef}
 \gamma_x = \lim_{m\to\infty} \beta_x(\kk_m).
\end{equation}
This limit exists as the sequence is increasing and bounded by $1$.
By monotone convergence,
\[
 \int_\sss \kk(x,y)\gamma_y\dd\mu(y) = \lim_{m\to\infty} \int_\sss\kk_m(x,y)\beta_y(\kk_m)\dd\mu(y)
\]
for every $x$, so
\begin{eqnarray}\nonumber
 \pgem\left(\int_\sss \kk(x,y)\gamma_y \dd\mu(y)\right)
&=&\lim_{m\to\infty}  \pgem\left(\int_\sss \kk_m(x,y)\beta_y(\kk_m) \dd\mu(y)\right) \\
&=& \lim_{m\to\infty} \beta_x(\kk_m) = \gamma_x.\label{gfe}
\end{eqnarray}
In other words, $\gamma_x$ also satisfies the functional
equation \eqref{fe}. Also, again by monotone convergence,
\begin{equation}\label{gplus}
 \lim_{m\to\infty} \bp(\kk_m) = \pgek\left(\int_\sss \kk(x,y)\gamma_y \dd\mu(y)\right).
\end{equation}
In sufficiently nice cases, this can be used
to establish the continuity of $\kk$, by showing that the functional
equation \eqref{fe}
has only one strictly positive solution, and using a simple special
form of $\kk_1$ to ensure that $\lim_{m\to\infty}\bp(\kk_m)\ge \bp(\kk_1)>0$.

\subsection{Case study: a power-law or `scale-free' graph}

One of the most studied examples
of the general model $\gnkx$ is the following:
let $\sss=(0,1]$, with $\mu$ Lebesgue measure, let $\kk(x,y)=c/\sqrt{xy}$,
and set $\xn_i=i/n$. Then the edges
of $\gnkx$ are present independently, and the probability
that the edge $ij$ is present is $c/\sqrt{ij}$, or, if $c>\sqrt{2}$,
the minimum of this quantity and $1$.
This random graph $\gnkx$ is the `mean-field' version
of the `scale-free' random graph introduced by Barab\'asi and Albert
in \cite{BAsc} as a model of the world-wide web,
and studied in many papers, especially in the computer science and
statistical physics literature.

For this kernel, $\alpha_\ell(\kk)=\infty$ for all $\ell\ge 2$:
the `early' vertices (with $i/n$ small) have large degree, and, while the average
degree is bounded as $n\to\infty$, the average square degree is not.
If $\ell$ is fixed, $i/n$ is bounded away from zero, and $\eps\to 0$ slowly,
then the expected number of paths from vertex $i$ to an `early' vertex
$j$ with $j\le \eps n$ is large (order $(\log n)^{\ell-1}$), and
one can check that \eqref{asum2} does not hold.
On the other hand, it is easy to check that the probability
that there is such a path tends to zero uniformly in $n$ as $\eps\to 0$.
Indeed, this probability is bounded by the expected number of {\em late--early}
paths of length at most $\ell$ starting at $i$, where a late--early path
is a path all of whose vertices have indices at least $\eps n$, apart
from the last vertex which does not; see Section 4 of \cite{Rsmall}.

Let $\kkeps$ be the `truncated' kernel that agrees with $\kk$ on $[\eps,1]^2$
and is zero otherwise; in the graph, this corresponds to deleting all
edges incident with early vertices. From the observation above,
if $\eps=\eps(n)\to 0$, then $o_p(n)$ vertices are within
graph distance $\ell$ of an early vertex,
so only $o_p(n)$ vertices of $\gnkx$ are incident
with a path of length $\ell$ in $\gnkx$ not present in $\gnxx{\kkeps}$.
Letting $\eps\to 0$ slowly enough, it is easy to approximate
the bounded kernels $\kkeps$ by regular finitary kernels (for example,
step functions), and so to deduce that \eqref{asum3} holds.
Thus
we shall be able to apply Theorem \ref{t3}, if we can establish
the required continuity of $\bp$.

Since we shall vary the parameter $c$, it will be convenient to write $\kk$ as $c\kk_0$,
where $\kk_0(x,y)=1/\sqrt{xy}$. The size of the giant
component in $\gnxx{c\kk_0}$ was found in \cite{Rsmall};
in Section 6.1 it was shown that the size is $\sigma(c)n+o_p(n)$,
where $\sigma(c)$ is the survival probability of the branching process $X_{c\kk_0}$.
(This is a special case of the main result of \cite{kernels}.)
In Section 6.2 of \cite{Rsmall}, this survival probability
is calculated in terms of the exponential integral. It turns out
that there is a giant component for {\em any} $c>0$, although
it is extremely small when $c$ is small: $\sigma(c)\sim 2e^{1-\gamma}\exp(-1/(2c))$
as $c\to 0$, where $\gamma$ is Euler's constant.

One might expect that
for $k$ fixed, there would be a $k$-core in $\gnxx{c\kk_0}$ for any $c>0$, perhaps
with size a constant fraction of that of the giant component when $c$ is small.
In fact, there is a positive threshold above which the $k$-core first appears.
Just above this threshold, the $k$-core is small: this is
in sharp contrast to $G(n,\la/n)$. In this result we write $\bp_k$
for $\bp$ when it is necessary to indicate the dependence on $k$.

\begin{theorem}\label{t4}
Let $\kk_0(x,y)=1/\sqrt{xy}$. For $c>0$,
let $\gnxx{c\kk_0}$ be the graph on $[n]=\{1,2,\ldots,n\}$ in which
edges are present independently, and the probability that $i$ and $j$
are joined is $\min\{c/\sqrt{ij},1\}$.
For each $k\ge 2$ we have
\[
 c_k\bb{\gnxx{c\kk_0}} = \bp_k(c\kk_0)n+o_p(n).
\]

If $k\ge 3$, then $\bp_k(c\kk_0)=0$ for $c\le (k-2)/2$, while
\begin{equation}\label{above}
 \bp_k(c\kk_0) \sim \frac{(k-1)!^{2/(k-2)}}{(k-1)(k-2)} \eps^{2/(k-2)}
\end{equation}
when $c=(1+\eps)(k-2)/2$ and $\eps\to 0$ from above.

If $k=2$, then $\bp_k(c\kk_0)>0$ for every $c>0$, and
\begin{equation}\label{2c}
 \bp_k(c\kk_0) \sim \frac{1}{2c}e^{2-2\gamma}\exp\left(-1/c\right)
\end{equation}
as $c\to 0$.
\end{theorem}

Thus, there is always a `giant' $2$-core; for small $c$ its size
is essentially the square of that of the giant component (times a factor $\Theta(1/c)$,
which is logarithmic in terms of the normalized size of the giant component).
For $k\ge 3$, the $k$-core emerges at a positive threshold, $c=(k-2)/2$,
and does so slowly. In the terminology of \cite{kernels}, for $k\ge 3$
the emergence of the $k$-core exhibits a phase transition of {\em exponent} $2/(k-2)$,
where this is the exponent of $\eps$ appearing above.
This contrasts with the exponent $0$ transition in $G(n,\la/n)$. (The
term `order' is often used in this context, but not always in the same way.)

\begin{proof}
As noted above, the graph $\gnxx{c\kk_0}$ satisfies \eqref{asum3}. Thus
the first statement follows from Theorem \ref{t3} once
we have established the required continuity of $\bp$ at $c\kk_0$. We
return to this later.

For the second and third statements,
let us write $\beta_x$ for $\beta_x(c\kk_0)$, which depends on $c$ and also on $k$.
Then equation \eqref{fe} becomes
\begin{equation}\label{fe2}
 \beta_x = \pgem\left(\int_0^1 \frac{c\beta_y}{\sqrt{xy}} \dd y\right).
\end{equation}
Define $A=A_k(c)$ by
\begin{equation}\label{Adef}
 A =  A_k(c) = \int_0^1 \frac{c\beta_y}{\sqrt{y}} \dd y.
\end{equation}
Then
\begin{equation}\label{B}
 \beta_x=\pgem(A/\sqrt{x}),
\end{equation}
so to determine $\beta_x$ it remains to determine $A$.
Substituting \eqref{B} into \eqref{Adef}, we see that $A=cf(A)$, where
\[
 f(B) =f_k(B) \= \int_0^1 \frac{\pgem(B/\sqrt{y})}{\sqrt{y}} \dd y.
\]
So far the argument is very similar to that in Section 6.2 of \cite{Rsmall};
indeed, for $k=2$ the event $\br$ is the event that the branching process survives,
so $\beta_x$ is the probability of survival starting
with a particle of type $x$; this is denoted $S_\infty(x)$
in \cite{Rsmall} (see page 919), where it is given by \eqref{B} with
\begin{equation}\label{a2c}
 A_2(c)\sim e^{1-\gamma}\exp(-1/(2c))
\end{equation}
as $c\to 0$, where $\gamma$ is Euler's constant.
It turns out that the case $k\ge 3$ is
much easier to handle.

Suppose that $k\ge 3$. Then, substituting $y=B^2x^{-2}$, so $x=B/\sqrt{y}$, we have
\[
 f_k(B) = \int_B^\infty \frac{x}{B} \pgem(x) 2B^2 x^{-3} \dd x  = 2Bg_k(B),
\]
where
\begin{equation}\label{gform}
 g_k(B) \= \int_B^\infty \pgem(x) x^{-2} \dd x.
\end{equation}
For any given $k\ge 3$, it is straightforward to calculate $f_k(B)$ explicitly.
Indeed,
$\pgem(x)=\sum_{t\ge k-1} x^t\exp(-x)/t!$,
and for $t\ge 2$ we have $\int_0^\infty x^{t-2}\exp(-x)\dd x=(t-2)!$.
Thus, for $k\ge 3$,
\[
 g_k(0) = \sum_{t\ge k-1} \frac{(t-2)!}{t!} = \sum_{t\ge k-1}\frac{1}{t(t-1)} = \frac{1}{k-2}.
\]
In particular, $g_3(0)=1$. It is easy to verify by differentiating that
\[
 g_3(B) = \frac{1-e^{-B}}B,
\]
so
\[
 f_3(B)= 2(1-e^{-B}) = 2B-B^2+O(B^3).
\]
Also, for each $k\ge 4$,
\[
 g_k(B)-g_{k+1}(B) = \int_B^\infty \frac{x^{k-3}}{(k-1)!}\exp(-x) \dd x,
\]
which is $\exp(-B)$ times a polynomial in $B$ that may be easily evaluated. Rather
than do this, let us note that, 
as $B\to 0$,
\[
 g_k(0)-g_k(B) = \int_{0}^B \frac{x^{k-1}}{(k-1)!}x^{-2} +O(x^{k-2}) \dd x
  = \frac{B^{k-2}}{(k-2)(k-1)!} + O(B^{k-1}),
\]
so
\begin{equation}\label{fasymp}
 f_k(B) = \frac{2B}{k-2}\left(1 - \frac{B^{k-2}}{(k-1)!} + O(B^{k-1})\right)  .
\end{equation}

Note that $f_k(B)/B=2g_k(B)$ is decreasing (from the form of \eqref{gform}).
Thus the equation $B=cf_k(B)$ has a positive solution if and only if
$cf_k'(0)>1$, i.e., if and only if $c>(k-2)/2$.
Furthermore, when $c>(k-2)/2$,
the solution $A_k(c)$ is unique. 

Let $c_0=c_0(k)=(k-2)/2$, and let $c=(1+\eps)c_0$ with $\eps>0$.
From \eqref{fasymp} and the fact that $f_k(B)/B$ is decreasing in $B$,
it follows that $A=A_k(c)\to 0$ as $\eps\to 0$. Furthermore,
as $\eps\to 0$,
\begin{eqnarray*}
 A = cf_k(A) &=& \frac{c}{c_0}A\left(1-(1+o(1))\frac{A^{k-2}}{(k-1)!}\right)\\
 &=& A + \eps A - (1+o(1))A\frac{A^{k-2}}{(k-1)!},
\end{eqnarray*}
so
\begin{equation}\label{a+}
 A = A_k(c) =A_k((1+\eps)c_0) \sim \left((k-1)!\eps\right)^{1/(k-2)}.
\end{equation}
 
Recalling that any solution to the functional equation \eqref{fe2} has the form
\eqref{B} with $A$ satisfying $A=cf(A)$, and that the probability
$\beta_x=\Pr_{c\kk,x}(\br)$ is given by the maximum solution to this equation,
we have shown that if $c\le c_0(k)$, then $\beta_x=0$ for all $x$,
while if $c> c_0(k)$, then
\[
 \beta_x= \pgem(A_k(c)/\sqrt{x}).
\]
Using \eqref{bp1} and \eqref{bp2}, we have $\bp(c\kk_0)=0$ if $c\le c_0(k)$.
Otherwise,
\begin{equation}\label{c+}
 \bp(c\kk_0) = \int_{x=0}^1\pgek(A_k(c)/\sqrt{x})\dd x.
\end{equation}
Let
\[
 h(B) = \int_{x=0}^1\pgek(B/\sqrt{x})\dd x = \int_{y=B}^\infty \pgek(y) 2B^2y^{-3} \dd y.
\]
As $B\to 0$ we have
\[
 \frac{h(B)}{2B^2} \to \int_{0}^\infty \sum_{t\ge k}\frac{y^{t-3}}{t!}e^{-y} \dd y
 = \sum_{t\ge k} \frac{(t-3)!}{t!} = \frac{1}{2(k-1)(k-2)}.
\]
Hence, from \eqref{a+} and \eqref{c+},
\[
 \bp\bb{(1+\eps)c_0\kk_0} \sim \frac{(k-1)!^{2/(k-2)}}{(k-1)(k-2)} \eps^{2/(k-2)},
\]
which is exactly \eqref{above}.

Returning to the case $k=2$, in this case, by reducing
to the exponential integral as in \cite{Rsmall}, one can show that
$h(B)\sim \log(1/B)B^2$ as $B\to 0$. Using \eqref{a2c}, \eqref{2c} follows.

To complete the proof of Theorem \ref{t4}, it remains to show that
there is a sequence $\kk_m$ of finite-type kernels
with $\sup_{x,y} \kk_m/(c\kk_0)<1$ for every $m$,
such that $\bp(\kk_m)\to\bp(c\kk_0)$.
We may assume that $c>(k-2)/2$, as otherwise
$\bp(c\kk_0)=0$ and there is nothing to prove.
Given a sequence $\kk_m$ tending up to $\kk$,
let $\gamma_x$ be defined by \eqref{gdef}. Then, from \eqref{gfe}, the function $x\mapsto\gamma_x$
also satisfies the functional equation \eqref{fe2}. We have shown above that this
functional equation has exactly two solutions, the zero function, and $\beta_x$.
Hence, writing $\gp$ for $\lim_{m\to\infty}\bp(\kk_m)$, from \eqref{gplus}
we have $\gp=0$ or $\gp=\bp(c\kk_0)$. We shall rule out the former case
by constructing the first element $\kk_1$ of our approximating sequence
suitably.
As a first step, we consider a bounded kernel that approximates $c\kk_0$.

Let $(k-2)/2<c'<c$ be fixed. Given $\eps>0$, let $\kkeps$ be the kernel given
by $\kkeps(x,y)=c'\kk_0(x,y)=c'/\sqrt{xy}$ if $x$, $y\ge \eps$, and $\kkeps(x,y)=0$
otherwise. Arguing exactly as for $c\kk_0$, we have
\[
 \beta_x(\kkeps) = \pgem(\Aeps/\sqrt{x})
\]
for $\eps\le x\le 1$, and $\beta_x(\kkeps)=0$ otherwise, where $\Aeps$
is the largest solution to $\Aeps=c'\feps(\Aeps)$, with
\[
 \feps(B) = \int_{\eps}^1 \frac{\pgem(B/\sqrt{y})}{\sqrt{y}} \dd y.
\]
Since $c'>(k-2)/2$, there is a $B>0$ with $c'f(B)>B$; fix such a $B$.
As $\eps\to 0$, we have $\feps(B)\upto f(B)$, so there is an $\eps>0$
with $c'\feps(B)>B$. It follows that $\Aeps>0$, and hence
that $\bp(\kkeps)>0$.
To complete the proof, let $\kk_1$ be a regular finitary kernel
with $c'\kkeps(x,y)\le \kk_1(x,y)\le \frac{1}{2}(c+c')\kk(x,y)$ for all $x,y$. Such
a kernel is easy to construct as $\kk$ is continuous and bounded
away from $0$ on the compact set $[\eps,1]^2$, where it coincides with $\kkeps$.
It is easy to construct a sequence of finite-type kernels $\kk_m$
with $\sup_{x,y}\kk_m/(c\kk_0)<1$, with $\kk_m$
tending up to $c\kk_0$ and starting with this particular $\kk_1$.
Then $\gp\ge \bp(\kk_1) \ge \bp(c'\kkeps)>0$. 
Since $\gp=0$ or $\gp=\beta(c\kk_0)$, we have $\gp=\beta(c\kk_0)$ as required.
\end{proof}

There is a sense in which Theorem \ref{t4} does not illustrate the
full power of Theorem \ref{t3}. Indeed, the kernel $c/\sqrt{xy}$
has a special property: it may be written as $\kk(x,y)=\phi(x)\phi(y)$
for some function $\phi$ on $\sss$.
Such kernels are called {\em rank 1} in \cite{kernels}.
The branching process corresponding to a rank 1 kernel is much simpler
than the general case: roughly speaking,
while the distribution of the {\em number} of children
of a particle depends on its type, the distribution of their {\em types} does not.

There is a minor variant of the model $\gnkx$, where the edge
probabilities are taken as $\pij=\kk(\xn_i,\xn_j)/(n+\kk(\xn_i,\xn_j))$,
so $\pij/(1-\pij)=\kk(\xn_i,\xn_j)/n$. As in \cite{kernels}, all our
results here apply equally to this variant. As noted in Section 16.4 of \cite{kernels},
using this variant, if $\kk$ has rank 1 then,
conditional on the degree sequence, $\gnkx$ is equally
likely to be any graph with the given degree sequence. Thus
the structure of $\gnkx$ is simpler than in the general case.
Also, $\gnkx$ is then closely related to random graphs
defined by first fixing a degree sequence (perhaps exactly,
or perhaps asymptotically), and then choosing a random graph with this
degree sequence.

Graphs of this form have been studied by many authors,
including Molloy and Reed \cite{MR1} in the general case, and Aiello,
Chung and Lu \cite{AielloCL} in the power-law case. The $k$-core
of such graphs has been studied by Janson and Luczak \cite{JL}
and by Fernholz and Ramachandran \cite{FR}; it is probable
that Theorem \ref{t4} could be proved by the methods of either
of these papers, although the calculation of $\bp(c\kk_0)$
must still be carried out.

In fact, Fernholz and Ramachandran studied the $k$-core in random graphs with
a given power-law degree sequence with exponent
$\alpha$: they assume a limiting fraction $i^{-\alpha}/\zeta(\alpha)$
of vertices with degree $i$ for each $i\ge 1$. When $\alpha=3$,
this model is very close to $\gnxx{c\kk_0}$, where the
degree distribution satisfies $\Pr(\mathrm{deg}(v)=i) \sim a i^{-3}$
as $i\to\infty$, for a constant $a$ depending on $c$.
Note, however, that, having fixed the degree exponent, the model considered here is much 
more flexible, due to the presence of the parameter $c$, which allows control
of the overall number of edges. For this reason,
the result in \cite{FR}, that \whp\ there is a $k$-core if $\alpha<3$,
and \whp\ there is no $k$-core when $\alpha\ge 3$,
gives no insight into the transition studied in Theorem \ref{t4}.
Roughly speaking, the $\alpha=3$ case of this result corresponds to showing that there
is no $3$-core in $\gnxx{c\kk_0}$ for a specific $c$. However,
the correspondence is not direct: in the model of \cite{FR}, the
presence or absence of the $k$-core depends very much on the entire degree
distribution, not just its asymptotics, and in $\gnxx{c\kk_0}$,
the distribution of the small degrees does not follow exactly a power-law.

It would be interesting to use Theorem \ref{t3} to compute the size of the $k$-core
in examples of $\gnkx$ where $\kk$ does not have rank one. A particularly
interesting case is the kernel $\kk(x,y)=c/\max\{x,y\}$ on $(0,1]^2$,
corresponding to the `uniformly grown random graph' proposed by Dubins
(as an infinite random graph) in 1984 (see~\cite{KW}). A closely related
model was introduced by 
Callaway, Hopcroft, Kleinberg, Newman and Strogatz~\cite{CHKNS} in 2001.
The emergence of the giant component in this model shows particularly interesting
behaviour: at $c=1/4$ there is an `infinite order' phase transition;
see~\cite{DMS-anomalous,Durrett,BJR,Rsmall}.
The functional equation \eqref{fe} is likely to be harder to handle in this case,
but may well still be tractable. Indeed, although the exact solution
could not be calculated, good bounds on the size of the giant component
just above the transition were obtained in \cite{Rsmall} by bounding the solution to
a related functional equation for the giant component in a
generalization of this model.

\medskip\noindent
{\bf Acknowledgements.} The author would like to thank Alan Frieze for suggesting
studying the $k$-core of $\gnkx$, and B\'ela Bollob\'as and Svante Janson
for helpful discussions on this topic.

\end{document}